\documentclass[10pt, a4paper]{article}

\usepackage{fullpage}
\usepackage{amssymb}
\usepackage{amsmath}
\usepackage{amsfonts}
\usepackage{subcaption}
\usepackage{multirow}
\usepackage{tikz}

\newtheorem{remark}{Remark}

\usepackage{xcolor}
\usepackage[colorlinks]{hyperref}
\definecolor{darkgreen}{rgb}{0.0,0.4,0.0}
\definecolor{darkred}{rgb}{0.6,0.0,0.0}
\definecolor{darkblue}{rgb}{0.0,0.0,0.5}
\definecolor{gray}{rgb}{0.5,0.5,0.5}
\definecolor{cyan}{rgb}{0.0,1.0,1.0}
\definecolor{darkcyan}{rgb}{0.0,0.5,0.5}
\definecolor{darkorange}{rgb}{0.8,0.4,0.0}
\definecolor{darkmargenta}{rgb}{0.5,0.0,0.5}
\definecolor{black}{rgb}{0.0,0.0,0.0}

\hypersetup{citecolor=blue, linkcolor=red, anchorcolor=darkblue,
filecolor=darkblue, urlcolor=darkblue}
\hypersetup{pdfauthor={Nicolae Suciu},pdftitle={Suciu}}

\begin{document}




\title{Global random walk solvers for fully coupled flow and transport in saturated/unsaturated porous media (extended version)}

\author{Nicolae Suciu$^{1, 2}\footnote{Corresponding author. \textit{Email adress: suciu@math.fau.de}}$, Davide Illiano$^3$, Alexander Prechtel$^{1}$, Florin A. Radu$^3$ \\ \\
\smallskip 
\textit{$^{1}$ Mathematics Department, Friedrich-Alexander University of Erlangen-N\"urnberg,}\\
\textit{Cauerstra{\ss}e. 11, 91058 Erlangen, Germany}\\
\textit{$^{2}$Tiberiu Popoviciu Institute of Numerical Analysis, Romanian Academy,}\\
\textit{Fantanele 57, 400320 Cluj-Napoca, Romania}\\
\textit{$^{3}$ Department of Mathematics, University of Bergen,}\\
\textit{ All\'{e}gaten 41, 5007 Bergen, Norway}}

\date{}
\maketitle

\begin{abstract}
In this article, we present new random walk methods to solve flow and transport problems in unsaturated/saturated porous media, including coupled flow and transport processes in soils, heterogeneous systems modeled through random hydraulic conductivity and recharge fields, processes at the field and regional scales. The numerical schemes are based on global random walk algorithms (GRW) which approximate the solution by moving large numbers of computational particles on regular lattices according to specific random walk rules. To cope with the nonlinearity and the degeneracy of the Richards equation and of the coupled system, we implemented the GRW algorithms by employing linearization techniques similar to the $L$-scheme developed in finite element/volume approaches. The resulting GRW $L$-schemes converge with the number of iterations and provide numerical solutions that are first-order accurate in time and second-order in space. A remarkable property of the flow and transport GRW solutions is that they are practically free of numerical diffusion. The GRW solutions are validated by comparisons with mixed finite element and finite volume solutions in one- and two-dimensional benchmark problems. They include Richards' equation fully coupled with the advection-diffusion-reaction equation and capture the transition from unsaturated to saturated flow regimes. For completeness, we also consider decoupled flow and transport model problems for saturated aquifers.
\end{abstract}

\begin{flushleft}

\end{flushleft}

\section{Introduction}
\label{intro}

The accuracy and the robustness of the numerical schemes is the primary requirement for reliable and meaningful results of the current efforts to improve the understanding of the complexity and interdependence of the flow and transport processes in subsurface hydrology through numerical investigations.
Numerical solvers for partial differential equations modeling individual or coupled processes are often used as basic elements in the formulation of the more complex problems of practical interest, such as parameter identification \cite{Franssenetal2009}, hydraulic tomography \cite{Bellinetal2020}, Monte Carlo approaches for systems with randomly distributed parameters \cite{Pasettoetal2011}, or upscaling for mutiphase flows in heterogeneous subsurface formations \cite{Hajibeygietal2020}, among others.

A central issue in subsurface hydrology is the need of robust and computationally efficient numerical models for partially saturated soil-groundwater systems. The transition between unsaturated and saturated zones is particularly challenging. In unsaturated flows the water content $\theta$ and the hydraulic conductivity $K$ depend nonlinearly on the pressure head $\psi$ through material laws based on experiments, as far as $\psi<0$. The evolution of $\psi$ is governed by the parabolic Richards' equation which degenerates to a (generally) linear elliptic equation (i.e. the equation for steady-state flow in aquifers) if $\psi\ge 0$ \cite{AltLuckhaus}. Since the regions where degeneracy takes place depend on the evolution of the pressure $\psi$ in time and space, they are not known {\it a priori}. To cope with the nonlinearity and degeneracy of the Richards' equation, different linearization methods are needed, such as the Newton scheme \cite{Schneid2000,Hajibeygietal2020,KnabnerandAngermann2003}, which is second-order convergent but converges only locally (requires a starting point close enough to the solution) or the more robust but only first-order convergent Picard, modified Picard, or $L$ schemes \cite{Slodicka,Pop2004,ListandRadu2016,Raduetal2017}.

Explicit and implicit schemes proposed for nonlinear flows in unsaturated regime provide solutions with comparable accuracy but are generally ambiguous to compare in terms of computing time. Since they do not need to solve systems of linear algebraic equations at every time step, explicit schemes are in principle faster \cite{Liuetal2020} but their speed may be seriously affected by the need to use very small time steps \cite{Haverkampetal1977,Caviedesetal2013,Alecsaetal2020}. The time step in explicit schemes is constraint by stability conditions \cite{Strikwerda2004,Liuetal2020,Caviedesetal2013} and has to be significantly reduced to ensure small local P\'{e}clet number (P\'{e}), defined with respect to the space step. Large (global) P\'{e} characterizes advection-dominated transport problems \cite{BauseandKnabner2004,Kuzmin2009}. In such cases, reducing the local P\'{e} is a remedy to avoid the numerical diffusion and the oscillatory behavior of the solution \cite{Raduetal2011}. The criterion of small local P\'{e} is also recommended for numerical schemes solving the pressure equation in saturated flows \cite{Gotovacetal2009} and, since Richards' equation has the structure of the advection-diffusion equation, the recommendation holds for the unsaturated flows as well.

Well known approaches to avoid the numerical diffusion are the particle tracking in continuous space and the discrete random walk on lattices \cite{Suciu2019}. The accuracy of these schemes is determined by the number of computational particles undergoing random jumps in continuous space or on discrete lattices. In random walk schemes, the increase of the computation time with the number of particles is simply avoided by randomly distributing the particles along the spatial directions with a global procedure, according to appropriate jump probabilities. In this way, one obtains a global random walk (GRW) which performs the spreading of all the particles from a given site with computational costs that are practically the same as for generating the jump of a single random walker in sequential procedures \cite{Vamosetal2003}. In particular cases (e.g., when using biased jump probabilities to account for variable coefficients or for advective displacements) the GRW algorithms are equivalent to explicit finite difference schemes with time step size constrained by stability requirements.  In unbiased GRW schemes for transport problems with variable coefficients, which still satisfy stability conditions, no restrictions on the time step are needed to reduce the local P\'{e} number, which renders the approach particularly efficient in large scale simulations of transport in groundwater (see \cite{Suciu2019} for details and examples).

The elliptic and parabolic equations governing the pressure head for flows in unsaturated/saturated porous media are essentially diffusion equations with second order operator in Stratonovich form. They can be recast as Fokker-Planck equations, with drift augmented by the row derivative of the coefficient tensor, and further solved by random walk approaches \cite{Suciu2019}. An alternative approach starts with a staggered finite difference scheme, further used to derive biased random walk rules governing the movement on a regular lattice of a system of computational particles. The particle density at lattice sites provides a numerical approximation of the pressure head solution. This approach has been already illustrated for flows in saturated porous media with heterogeneous hydraulic conductivity \cite{Alecsaetal2020,Suciu2020}.

In this article, we present new GRW schemes for nonlinear and non-steady flows in soils which model the transition from unsaturated to saturated regime in a way consistent with the continuity of the constitutive relationships $\theta(\psi)$ and $K(\psi)$. Following \cite{ListandRadu2016,Raduetal2017}, the nonlinearity of the Richards equation is solved with an iterative procedure similar to the $L$-scheme used in finite element/volume approaches. Numerical tests demonstrate the convergence of the $L$-scheme for unsaturated/saturated flows. For fully saturated flow regime with constant water content $\theta$ and time independent boundary conditions the GRW $L$-scheme is equivalent to a transient finite difference scheme.

Coupled flow and reactive transport problems for partially saturated soils rise new stability and consistency issues and demand augmented computational resources. Our GRW approach in this case consists of coupling the flow solver described above with existing GRW transport solvers \cite{Suciu2019} adapted for nonlinear problems, which are implemented as $L$-schemes as well. The flow and transport solvers are coupled via an alternating splitting procedure \cite{Illiano2020} which successively iterates the corresponding $L$-schemes until the convergence of the pressure head and concentration solutions is reached, within the same tolerance, at every time step. Code verification tests using analytical manufactured solutions are employed to verify the convergence of the iterations and the accuracy of the splitting scheme.

 The GRW scheme for one-dimensional solutions of the Richards equation, which captures the transition from unsaturated to saturated flow regimes is validated by comparisons with solutions provided by {\sc Richy} software, based on the mixed finite element method (MFEM), with backward Euler discretization in time and Newton linearization, developed at the Mathematics Department of the Friedrich-Alexander University of Erlangen-N\"urnberg \cite{Schneidetal2000,Schneid2000}. For the particular case of unsaturated flows, the one-dimensional GRW solutions are also tested by comparisons with experimental data \cite{Zadeh2011,Zambraetal2012} and exact solutions published in the literature \cite{Warricketal1985,Watsonetal1995}. The two-dimensional GRW solutions are compared on benchmark problems with the results obtained thanks to two-point flux approximation (TPFA) finite volume solvers using backward Euler discretization in time and $L$-scheme linearization \cite{Illiano2020}. The TPFA solvers are implemented in MRST, the MATLAB Reservoir Simulation Toolbox \cite{mrst}.

The paper is organized as follows. Section~\ref{1Dalg} presents the GRW algorithm and the linearization approach for one-dimensional flow problems. The one-dimensional solver is further validated through comparisons with MFEM solutions, experimental data, and exact analytical solutions in Section~\ref{1Dvalid}. Two dimensional GRW algorithms for fully coupled and decoupled flow and transport problems are introduced in Section~\ref{2D}. Code verification tests and validation via comparisons with TPFA solutions for benchmark problems are presented in Section~\ref{2Dvalid}. Some examples of flow and transport solutions for groundwater problems at the field and regional scale are presented in Section~\ref{2Dreg_scale}. The main conclusions of this work are finally presented in Section~\ref{concl}. Appendix \ref{appA} contains estimations of computational order of convergence of the GRW $L$-scheme. GRW codes implemented in Matlab for model problems considered in this article are stored in the Git repository \href{https://github.com/PMFlow/RichardsEquation}{RichardsEquation} \cite{Suciuetal2021}.

\section{One-dimensional GRW algorithm for unsaturated/saturated flow in soils}
\label{1Dalg}

We consider the water flow in unsaturated/saturated porous media described by the one-dimensional Richards equation \cite{Haverkampetal1977,Schneidetal2000,KnabnerandAngermann2003} in the space-time domain $[0,L_z]\times[0,T]$,
\begin{equation}\label{eq1}
\frac{\partial \theta(\psi)}{\partial t}-\frac{\partial}{\partial z}\left[K(\theta(\psi))\frac{\partial}{\partial z}(\psi+z)\right]=0,
\end{equation}
where $\psi(z,t)$ is the pressure head expressed in length units, $\theta$ is the volumetric water content, $K$ stands for the hydraulic conductivity of the medium, and $z$ is the height oriented positively upward. According to (\ref{eq1}), the water flux given by Darcy's law is $q=-K(\theta(\psi))\frac{\partial}{\partial z}(\psi+z)$.

To design a GRW algorithm, we start with the staggered finite difference scheme with backward discretization in time which approximates the solution of Eq.~(\ref{eq1}) at positions $z=i\Delta z$, $i=1,\ldots, I$, $I=L_z/\Delta z$, and time points $t=k\Delta t$, $k=1,\ldots, T/\Delta t$, according to
\begin{align}\label{eq4}
\theta(\psi_{i,k})-\theta(\psi_{i,k-1})=
\frac{\Delta t}{{\Delta z}^2}&\left\{\left[K(\psi_{i+1/2,k})(\psi_{i+1,k}-\psi_{i,k})-
K(\psi_{i-1/2,k})(\psi_{i,k}-\psi_{i-1,k})\right]\right.\nonumber\\
&\left.+\left(K(\psi_{i+1/2,k})-K(\psi_{i-1/2,k})\right)\Delta z\right\}.
\end{align}

To cope with the double nonlinearity due to the dependencies $K(\theta)$ and $\theta(\psi)$ we propose an explicit scheme similar to the linearization approach known as ``$L$-scheme'', originally developed for implicit methods \cite{Pop2004,ListandRadu2016,Raduetal2017}. The approach consists of the addition of a stabilization term $L(\psi_{i,k}^{s+1}-\psi_{i,k}^{s})$, $L=const$, in the left-hand side of (\ref{eq4}) and of performing successive iterations $s=1,2,\ldots $ of the modified scheme until the discrete $L^2$ norm of the solution $\psi_{k}^{s}=(\psi_{i,k}^{s}, \ldots, \psi_{I,k}^{s})$ verifies
\begin{equation}\label{eq5}
\|\psi_{k}^{s}-\psi_{k}^{s-1}\|\leq \varepsilon_a + \varepsilon_r\|\psi_{k}^{s}\|
\end{equation}
for some given tolerances $\varepsilon_a$ and $\varepsilon_r$. The adapted $L$-scheme reads
\begin{align}\label{eq6}
\psi_{i,k}^{s+1}=& \left[1-(r_{i+1/2,k}^{s}+r_{i-1/2,k}^{s})\right]\psi_{i,k}^{s} + r_{i+1/2,k}^{s}\psi_{i+1,k}^{s} + r_{i-1/2,k}^{s}\psi_{i-1,k}^{s} \nonumber\\
&+ \left(r_{i+1/2,k}^{s}-r_{i-1/2,k}^{s}\right)\Delta z
-\left(\theta(\psi_{i,k}^{s})-\theta(\psi_{i,k-1})\right)/L,
\end{align}
where
\begin{equation}\label{eq7}
r_{i\pm 1/2,k}^{s}=K(\psi_{i\pm 1/2,k}^{s})\Delta t/(L{\Delta z}^2).
\end{equation}
For fixed time step $k$, the iterations start with the solution after the last iteration at the previous time $k-1$, $\psi_{i,k}^{1}=\psi_{i,k-1}$, $i=1,\ldots, I$. Note that, unlike implicit $L$-schemes (e.g., \cite{Slodicka,Pop2004,ListandRadu2016}), the explicit scheme (\ref{eq6}) uses forward increments of $\psi$. In this way, the solution $\psi_{i,k}^{s+1}$ is obtained from values of $\psi$ and $r$ evaluated at the previous iteration, without solving systems of algebraic equations.

The solution $\psi_{i,,k}^{s}$ is further represented by the distribution of $\mathcal{N}$ computational particles at the sites of the one-dimensional lattice, $\psi_{i,k}^{s}\approx n_{i,k}^{s} a /\mathcal{N}$, with $a$ being a constant equal to a unit length, and the $L$-scheme (\ref{eq6}) becomes
\begin{equation}\label{eq8}
n_{i,k}^{s+1}= \left[1-\left(r_{i+1/2,k}^{s}+r_{i-1/2,k}^{s}\right)\right]n_{i,k}^{s} + r_{i+1/2,k}^{s}n_{i+1,k}^{s} + r_{i-1/2,k}^{s}n_{i-1,k}^{s} + \left\lfloor \mathcal{N}f^{s}\right\rfloor,
\end{equation}
where the source term is defined as $f^{s}=\left(r_{i+1/2,k}^{s}-r_{i-1/2,k}^{s}\right)\Delta z-\left[\theta(n_{i,k}^{s})-\theta(n_{i,k-1}\right]/L$ and $\left\lfloor\cdot\right\rfloor$ denotes the floor function.

The physical dimension of the parameter $L$ of the scheme is that of an inverse length unit to ensure that $r_{i\pm 1/2,k}^{s}$ defined by (\ref{eq7}) are non-dimensional parameters, as needed in random walk approaches \cite{Vamosetal2003,Suciu2019}. By imposing the constraint $r_{i\pm 1/2,k}^{s}\leq 1/2$, the parameters $r_{i\pm 1/2,k}^{s}$ can be thought of as biased jump probabilities. Hence, the contributions to $n_{i,k}^{s+1}$ from neighboring sites $i\pm 1$ summed up in (\ref{eq8}) can be obtained with the GRW algorithm which moves particles from sites $j$ to neighboring sites $i=j\mp 1$ according to the rule
\begin{eqnarray}\label{eq9}
n_{j,k}^{s}=\delta n_{j,j,k}^{s}+\delta n_{j-1,j,k}^{s}+\delta n_{j+1,j,k}^{s}.
\end{eqnarray}
For consistency with (\ref{eq8}), the quantities
$\delta n^{s}$ in (\ref{eq9}) have to satisfy in the mean \cite[Sect. 3.3.4.1]{Suciu2019},
\begin{equation}\label{eq10}
\overline{\delta n_{j,j,k}^{s}}=\left[1-\left(r_{j-1/2,k}^{s}+r_{j+1/2,k}^{s}\right)\right]\overline{n_{j,k}^{s}},\;\;
\overline{\delta n_{j\mp 1/2,j,k}^{s}}=r_{j\mp 1/2,k}^{s}\overline{n_{j,k}^{s}}.
\end{equation}

The quantities $\delta n^{s}$ are binomial random variables approximated by using the unaveraged relations (\ref{eq10}) for the mean, summing up the reminders of multiplication by $r$ and of the floor function $\lfloor \mathcal{N}f^s\Delta t\rfloor$, and allocating one particle to the lattice site where the sum reaches the unity.

\begin{remark}\label{remEM}
The finite difference $L$-scheme (\ref{eq6}) and the derived GRW relation (\ref{eq8}) are explicit iterative schemes for Richards equation in mixed form (\ref{eq1}). The essential difference of the $L$-schemes from explicit schemes in mixed formulation designed for unsaturated regime \cite{Haverkampetal1977,Caviedesetal2013,Liuetal2020} is that they apply to both unsaturated and saturated flow conditions.
\end{remark}

\begin{remark}\label{remFD} For fixed $k$ and iteration index $s$  interpreted as time, with $\Delta s=\Delta t$, the GRW relation (\ref{eq8}) is equivalent to a consistent forward-time central-space one-dimensional finite difference scheme for the following equation with variable coefficient $K$ and source term $f$,
\begin{equation}\label{eqFD}
\frac{\partial \psi}{\partial s}-\frac{\partial}{\partial z}\left[\frac{K(\theta(\psi))}{L}\frac{\partial \psi}{\partial z}\right]=f(\psi),
\end{equation}
with the right-hand side given by
\[
f(\psi)=\frac{1}{L}\left(\frac{\partial K(\theta(\psi(z,t,s)))}{\partial z}-\frac{\partial\theta(\psi(z,t,s))}{\partial t}\right).
\]
The equation is solved for $\psi(z,t,s)$ with initial condition $\psi(z,t,0)=\psi(z,t-\Delta t)$ and boundary conditions of the original problem for Eq. (\ref{eq1}). In the particular case of saturated regime with $\theta=const$, $f$ becomes independent of $\psi$ (possibly space-time variable through $\partial K/\partial z$) and the simplified scheme obtained from (\ref{eq8}) is easily seen to be first-order accurate in time and second-order in space. The jump probability $r$ corresponds to the von Neumann stability parameter. The constraint $r\leq 1/2$, fulfilled for every $s$ and space index $i$ (supplemented by Duhamel's superposition theorem if $f\neq 0$), ensures the stability of the scheme. Since the scheme is also consistent, the Lax-Richtmyer equivalence theorem implies its convergence \cite[Sects. 6.5, and 9.3]{Strikwerda2004}. Moreover, for smooth, or at least Lebesgue integrable initial data, the solutions converge with the order of accuracy of the scheme \cite[Sect. 10.4 ]{Strikwerda2004}. The convergence of the numerical solution provided by this scheme implies the convergence of the iterations, independently of the value of the stabilization parameter $L$.
\end{remark}

\begin{remark}\label{remFsat}
Consider again the saturated regime, $\theta=const$, with space-variable hydraulic conductivity $K$ and a given source term $f$. With the parameter $L$ set to $L=1/a$, after disregarding the time index $k$ the scheme (\ref{eq8}) solves the following equation for the hydraulic head $h=\psi+z$,
\begin{equation}\label{eqFsat}
\frac{1}{a}\frac{\partial h}{\partial s}-\frac{\partial}{\partial z}\left[K\frac{\partial h}{\partial z}\right]=f.
\end{equation}
For boundary conditions independent of $s$, the solution of Eq. (\ref{eqFsat}) approaches a steady-state regime corresponding to the saturated flow (see also \cite{Alecsaetal2020,Suciu2020}). The convergence of the scheme follows from Remark1~\ref{remFD}.
\end{remark}

\section{Validation of the one-dimensional GRW flow algorithm}
\label{1Dvalid}

\subsection{Transition from unsaturated to saturated flow regime}
\label{1Dvalid_degenerate}

The one-dimensional algorithm for flow in unsaturated/saturated soils is validated in the following by comparisons with MFEM solutions obtained with the {\sc Richy} software \cite{Schneidetal2000,Schneid2000}.
For this purpose, we solve one-dimensional model-problems for the vertical infiltration of the water through both homogeneous and non-homogeneous soil columns \cite{SrivastavaandYeh1991}, previously used in \cite{Schneid2000} to assess the accuracy and the convergence of the MFEM solutions.

We consider the domain $z\in[0,2]$ and the boundary conditions 
specified by a constant pressure $\psi(0,t)=\psi_0$ at the bottom of the soil column and a constant water flux $q_0$ at the top. Together, these constant conditions determine the initial pressure distribution $\psi(z,0)$ as solution of the steady-state flow problem. For $t>0$, the pressure $\psi_0$ is kept constant, at the bottom, and the water flux at the top of the column is increased linearly from $q_0$ to $q_1$ until $t\leq t_1$ and is kept constant for $t>t_1$.

\begin{figure}[h]\centering
\includegraphics[width=0.45\linewidth]{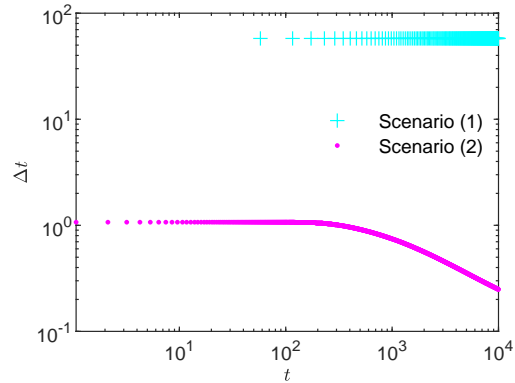}
\caption{\label{fig_tdt}Time steps for Scenario (1) and Scenario (2).}
\end{figure}

\begin{figure}[h]
\begin{minipage}[t]{0.45\linewidth}\centering
\includegraphics[width=\linewidth]{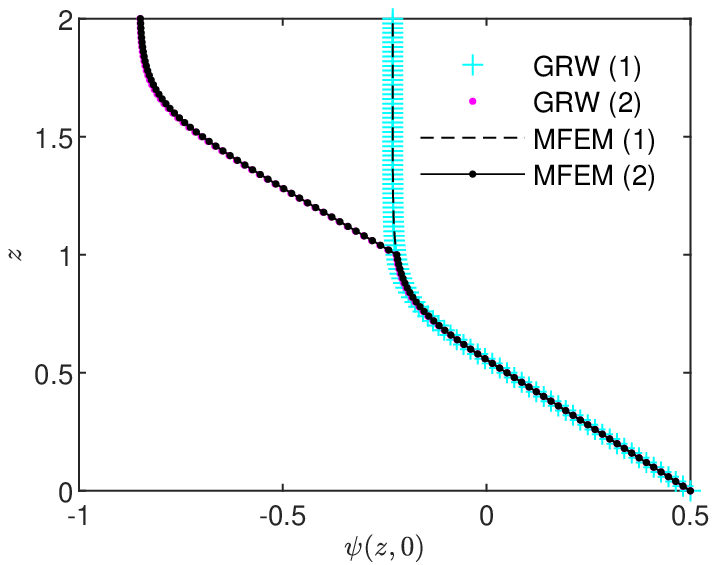}
\caption{\label{fig1}Initial condition for pressure head computed by GRW and MFEM codes.}
\end{minipage}
\hspace*{0.1in}
\begin{minipage}[t]{0.45\linewidth}\centering
\includegraphics[width=\linewidth]{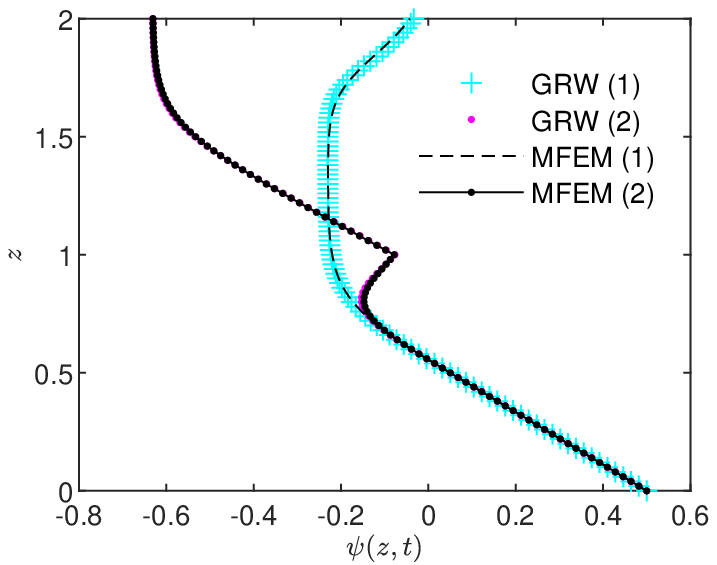}
\caption{\label{fig2}Pressure head solutions at $t=10^4$ seconds computed by GRW and MFEM codes.}
\end{minipage}
\end{figure}

\begin{figure}[h]
\begin{minipage}[t]{0.45\linewidth}\centering
\includegraphics[width=\linewidth]{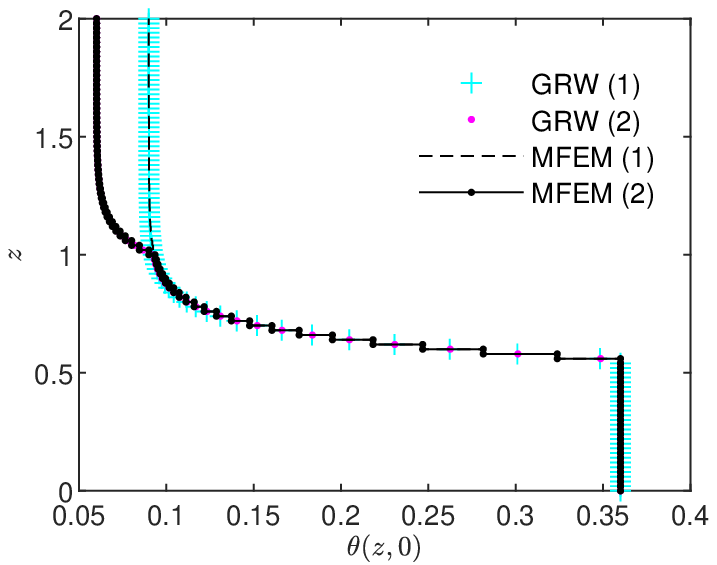}
\caption{\label{fig3}Initial water content computed by GRW and MFEM codes.}
\end{minipage}
\hspace*{0.1in}
\begin{minipage}[t]{0.45\linewidth}\centering
\includegraphics[width=\linewidth]{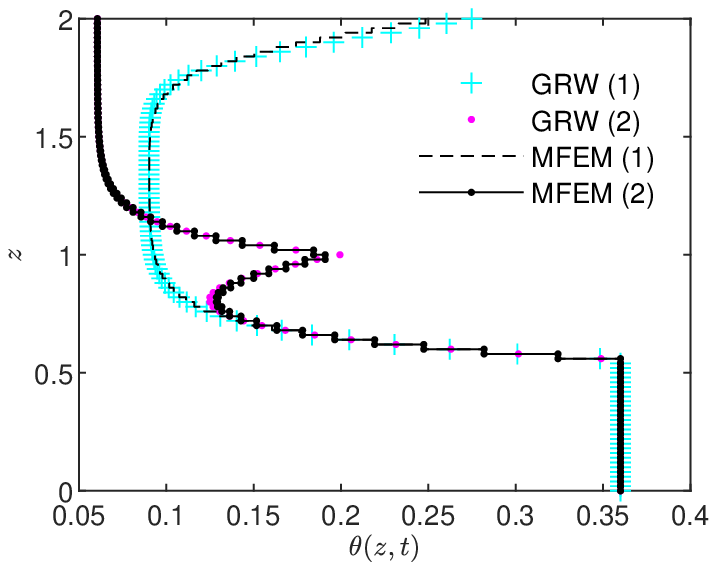}
\caption{\label{fig4}Water content solutions at $t=10^4$ seconds computed by GRW and MFEM codes.}
\end{minipage}
\end{figure}

\begin{figure}[h]
\begin{minipage}[t]{0.45\linewidth}\centering
\includegraphics[width=\linewidth]{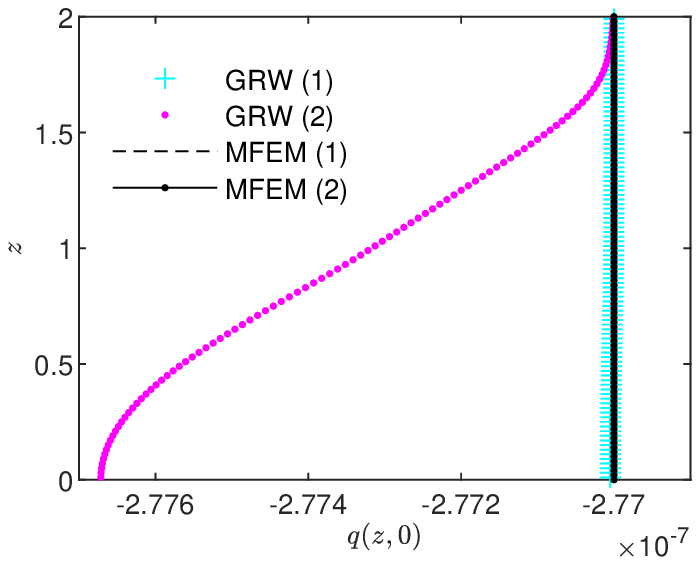}
\caption{\label{fig5}Initial water flux computed by GRW and MFEM codes.}
\end{minipage}
\hspace*{0.1in}
\begin{minipage}[t]{0.45\linewidth}\centering
\includegraphics[width=\linewidth]{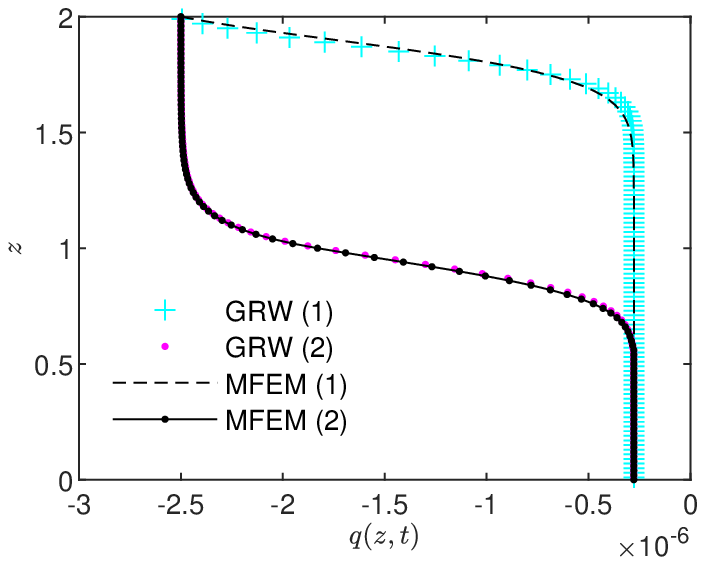}
\caption{\label{fig6}Water flux solutions at $t=10^4$ seconds computed by GRW and MFEM codes.}
\end{minipage}
\end{figure}

\begin{table}[h]
\centering
\caption{Error norms of the initial conditions.}
\label{table01}
\begin{tabular}{ c c c c c }
\hline
& $\varepsilon_{\psi_0}$ & $\varepsilon_{\theta_0}$ & $\varepsilon_{q_0}$ \\ \hline
Scenario (1) & 5.07e-03 & 1.85e-02 & 8.39e-06 \\ \hline
Scenario (2) & 4.55e-03 & 1.94e-02 & 1.46e-03 \\ \hline
\end{tabular}
\end{table}
\begin{table}[h]
\centering
\caption{Error norms of the GRW solutions.}
\label{table1}
\begin{tabular}{ c c c c c }
\hline
& $\varepsilon_{\psi}$ & $\varepsilon_{\theta}$ & $\varepsilon_{q}$ \\ \hline
Scenario (1) & 1.81e-02 & 2.20e-02 & 3.50e-02 \\ \hline
Scenario (2) & 5.20e-03 & 2.35e-02 & 2.07e-02 \\ \hline
\end{tabular}
\end{table}

For the unsaturated regions ($\psi< 0$) we consider the constitutive relationships given by the simple exponential model \cite{Gardner1958}
\begin{equation}\label{eq2}
\theta(\psi)=\theta_{res}+(\theta_{sat}-\theta_{res})e^{\alpha\psi},
\end{equation}
\begin{equation}\label{eq3}
K(\theta(\psi))=K_{sat}\frac{\theta(\psi)-\theta_{res}}{\theta_{sat}-\theta_{res}},
\end{equation}
where $\theta=\theta_{sat}$ and $K=K_{sat}$ denote the constant water content respectively the constant hydraulic conductivity in the saturated regions ($\psi\geq 0$) and $\theta_{res}$ is the residual water content. The more complex and physically sounded van Genuchten-Mualem parameterization model will be used in the next subsection for a comparison with measurements in unsaturated flow regime and for more complex two-dimensional problems in the following sections.

The flow problem for Eq.~(\ref{eq1}) with the parameterization (\ref{eq2}-\ref{eq3}) is solved in two Scenarios: (1) homogeneous soil, with $K_{sat} = 2.77\cdot 10^{-6}$, $\theta_{res}=0.06$, $\theta_{sat}=0.36$, $\alpha=10$, $q_0=2.77\cdot 10^{-7}$, $q_1=2.50\cdot 10^{-6}$, which are representative for a sandy soil, and (2) non-homogeneous soil, with the same parameters as in Scenario (1), except the saturated hydraulic conductivity, which takes two constant values, $K_{sat}= 2.77\cdot 10^{-6}$ for $z<1$ and $500 K_{sat}$ for $z\ge 1$ (modeling, for instance, a column filled with sand and gravel). To capture the transition from unsaturated to saturated regime, the pressure at the bottom boundary is fixed at $\psi_0=0.5$. For the parameters of the one-dimensional flow problems solved in this section we consider meters as length units and seconds as time units. The simulations are conducted up to $T=10^4$ (about 2.78 hours) and the intermediate time is taken as $t_1=T/10^2$.

We consider a uniform GRW lattice with $\Delta z=10^{-2}$, equal to the length of the linear elements in the MFEM solver. The GRW computations are initialized by multiplying the initial condition by $\mathcal{N}=10^{24}$ particles. Since, as shown by (\ref{eq3}), the hydraulic conductivity varies in time, the length of the time step determined by (\ref{eq7}) for the maximum of $K$ at every time iteration and by specifying a maximum $r_{max}=0.8$ of the parameter $r_{i\pm 1/2,k}$ may vary in time (see Fig.~\ref{fig_tdt}). The parameter of the regularization term in the $L$-scheme is set to $L=1$ for the computation of the initial condition (solution of the stationary problem, i.e. for $\partial\theta/\partial t=0$ in (\ref{eq1})) and to $L=2$ for the solution of the non-stationary problem. In both cases, the convergence criterion (\ref{eq5}) is verified by choosing $\varepsilon_a=0$ and a relative tolerance $\varepsilon_r=10^{-9}$.

The comparison with the MFEM solutions presented in Figs.~\ref{fig1}-\ref{fig6} shows a quite good accuracy of the GRW solutions for pressure, water content, and water flux. The relative errors, computed with the aid of the $L^2$ norms by $\varepsilon_{\psi}=\|\psi^{GRW}-\psi^{MFEM}\|/\|\psi^{MFEM}\|$, and similarly for $\theta$, $q$, and the solutions $\psi_0$, $\theta_0$ and $q_0$ of the steady state problem, are presented in Tables~\ref{table01}~and~\ref{table1}.

\begin{figure}[h]
\begin{minipage}[t]{0.45\linewidth}\centering
\includegraphics[width=\linewidth]{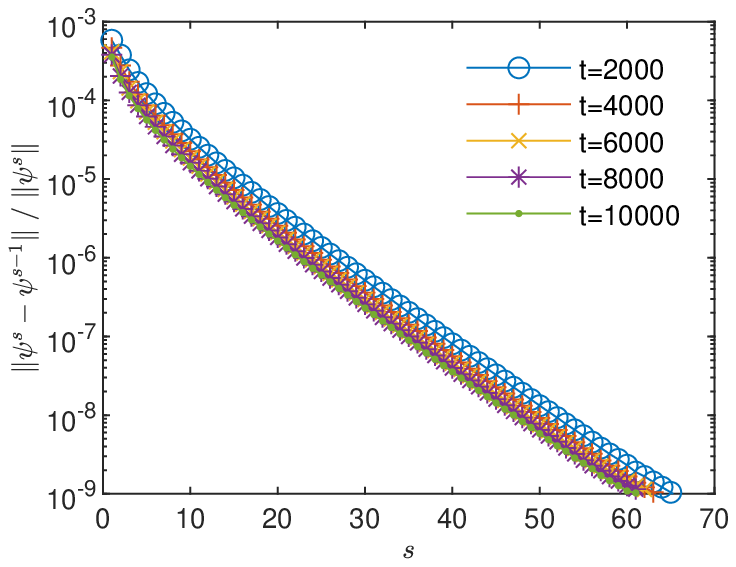}
\caption{\label{fig7}Convergence of the $L$-scheme implementation of the GRW solver in Scenario (1).}
\end{minipage}
\hspace*{0.1in}
\begin{minipage}[t]{0.45\linewidth}\centering
\includegraphics[width=\linewidth]{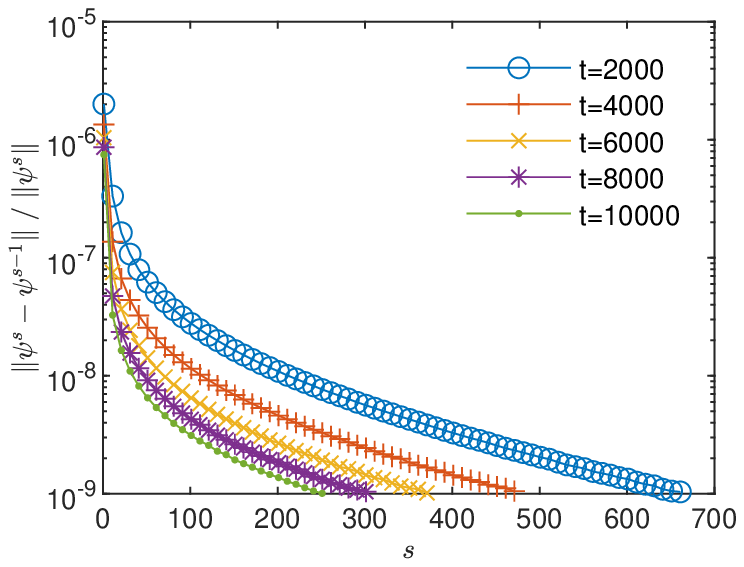}
\caption{\label{fig8}Convergence of the $L$-scheme implementation of the GRW solver in Scenario (2).}
\end{minipage}
\end{figure}

The $L$-scheme converges with speeds depending on the problem. To solve the problem for the initial condition, one needs $3.5\cdot 10^4$ iterations in Scenario (1) and $6.5\cdot 10^6$ iterations in Scenario (2). Instead, to solve the non-stationary problem for a final time $T=10^4$, one needs about 70 iterations in Scenario (1) and about 700 iterations in Scenario (2) (see Fig.~\ref{fig7} and Fig.~\ref{fig8}). The convergence of the iterative GRW $L$-scheme can be further investigated through assessments of the computational order of convergence of the sequence of successive correction norms $\|\psi_{k}^{s}-\psi_{k}^{s-1}\|$ \cite{Catinas2019,Catinas2020}. Estimations provided in Appendix~\ref{appA} indicate a linear convergence for Scenario (1) but only a power law convergence $\sim s^{-1}$, which is slower than the linear convergence \cite{Catinas2020}, for Scenario (2).

Supplementary tests done in Scenario 1 have shown that the convergence behavior depends on the value of the constant $L$. For $1\le L\le1.70$ the solution oscillates and the convergency indicator (\ref{eq5}) levels of at relatively large errors, between $10^{-1}$ and $10^{-3}$, then, the results improve with increasing $L$. For $1.96\le L\le 445.24$, the results are practically identical, with a plateau $\varepsilon_r\sim 10^{-10}$ of the relative errors given by (\ref{eq5}), reached after a number of iterations $s$ increasing from 70 to 10000. The upper bound of $L$ is the value of $L_\theta=\sup_{\psi}|\theta'(\psi)|$, computed with the parameters of the considered scenario. For finite element approaches it has been shown that $L\geq L_\theta$ is a sufficient condition to ensure the convergence of the $L$-scheme \cite{ListandRadu2016}. For the present random walk based scheme, it seems that $L_\theta$ does not play a special role. The tests show that increasing $L$ above the value which ensures the convergence of the $L$-scheme with a desired accuracy only results in increasing number of iterations and more computing time. The parameter $L$ has to be established experimentally by checking the convergence and, as highlighted by the examples presented in Section~\ref{2Dvalid} below, it depends on the complexity of the problem to be solved.

\subsection{Comparison with experiments and exact solutions for unsaturated flows}
\label{1Dvalid_unsaturated}

An experiment consisting of free drainage in a 600 cm deep lysimeter filled with a material with silty sand texture conducted at the Los Alamos National Laboratory \cite{Abeele1984} is often used to validate one-dimensional schemes for unsaturated flows (see e.g., \cite{Zadeh2011,Zambraetal2012,Caviedesetal2013}). This example is provided with the {\sc Hydrus 1D} software \cite{Simuneketal2008}, which is also used for validation purposes in the papers cited above.

The relationships defining the water content $\theta(\psi)$ and the hydraulic conductivity $K(\theta(\psi))$ are given by the van Genuchten-Mualem model
\begin{equation} \label{theta}
\Theta(\psi) = \begin{cases} \left(1+(-\alpha \psi)^n\right)^{-m}, &\psi < 0 \\
1, &\psi \geq 0,
\end{cases}
\end{equation}
\begin{equation} \label{K}
K(\Theta(\psi)) = \begin{cases} K_{sat} \Theta(\psi)^{\frac{1}{2}} \left[1-\left(1-\Theta(\psi)^\frac{1}{m}\right)^m \right]^2, &\psi < 0 \\
K_{sat}, &\psi \geq 0,
\end{cases}
\end{equation}
where $\theta_{res}$, $\theta_{sat}$, and $K_{sat}$ represent the same parameters as for the exponential model considered in Section~\ref{1Dvalid_degenerate}, $\Theta = (\theta - \theta_{res})/(\theta_{sat} - \theta_{res})$ is the normalized water content, and $\alpha$, $n$ and $m=1-1/n$ are model parameters depending on the soil type.

\begin{figure}[h]
\begin{minipage}[t]{0.45\linewidth}\centering
\includegraphics[width=\linewidth]{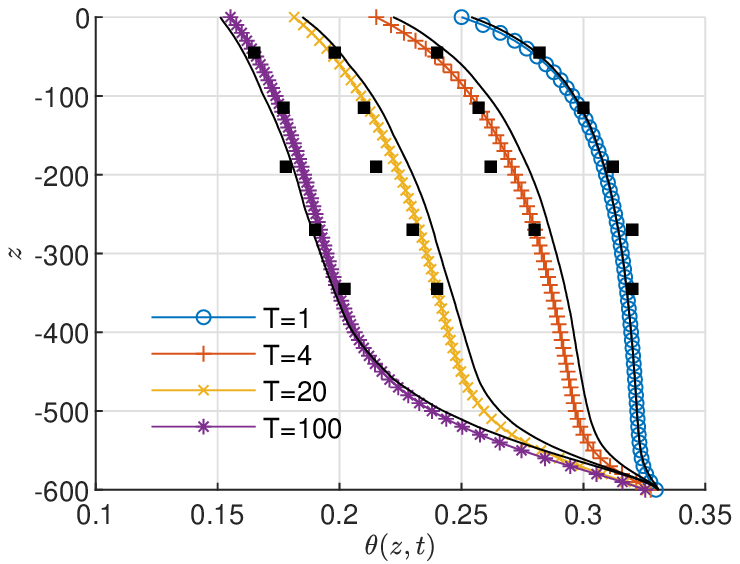}
\caption{\label{fig_drainage}Spatiotemporal distribution of the water content during the drainage experiment simulated by the GRW scheme. Continuous black lines represent the solution provided by {\sc Hydrus 1D} model. Black markers correspond to measurements picked-up from \cite[Fig. 2]{Zambraetal2012}.}
\end{minipage}
\hspace*{0.1in}
\begin{minipage}[t]{0.45\linewidth}\centering
\includegraphics[width=\linewidth]{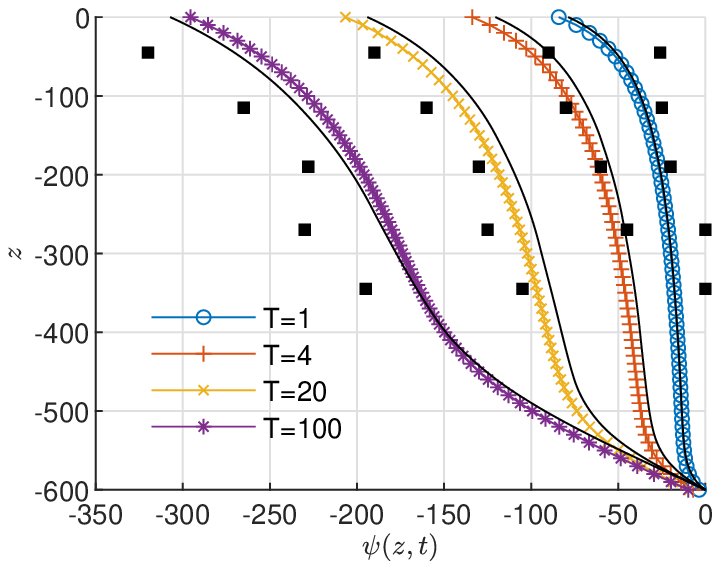}
\caption{\label{fig_drainage_p}Spatiotemporal distribution of the pressure head during the drainage experiment simulated by the GRW scheme. Continuous black lines represent the solution provided by {\sc Hydrus 1D} model. Black markers correspond to measurements picked-up from \cite[Fig. 4]{Zadeh2011}.}
\end{minipage}
\end{figure}

With the parameters given in the {\sc Hydrus 1D} example, $\theta_{res}=0.0$, $\theta_{sat}=0.331$, $K_{sat} = 25$ cm/d, $\alpha=0.0143$ cm$^{-1}$, $n=1.5$, for initial and boundary conditions for free drainage given by $\psi(z,0)=0$ cm and $q(0,t)=0$ cm/d \cite{Zadeh2011}, the solutions
 provided by the GRW $L$-scheme (\ref{eq8}-\ref{eq10}) for simulation times from 1 d to 100 d are obtained with stabilization parameter $L=0.5$ after a number of 9 to 35 iterations (tolerance specified by $\varepsilon_a=\varepsilon_r=5\cdot 10^{-6}$ in (\ref{eq5})). The spatial resolution is set to $\Delta z=10$ cm, while the time step varies slightly between $10^{-2}$ d and $3.16\cdot 10^{-2}$ d, according to (\ref{eq7}). The results are compared in Figs.~\ref{fig_drainage}~and~\ref{fig_drainage_p} with {\sc Hydrus 1D} results and experimental data. The water content profiles (Fig.~\ref{fig_drainage}) are quite close to measurements and similar to those presented in \cite{Zambraetal2012,Caviedesetal2013}. The pressure profiles  (Fig.~\ref{fig_drainage_p}) deviate from experiment, mainly for $T=1$ d and $T=100$ d, with approximately the same amount as in \cite[Fig. 12]{Caviedesetal2013}. An improved prediction of the pressure profiles is obtained in \cite{Zadeh2011} with slightly modified parameters of the van Genuchten-Mualem model, but with the price of larger deviations  for the water content.

The $\theta$-based form of Richards equation has shown significantly improved performance in numerical schemes for unsaturated flows in spatially homogeneous soils (e.g., constant $K_{sat}$), especially in modeling infiltration into dry media \cite{Zadeh2011}, and is well suited to analytical approaches \cite{Philip1969,Warricketal1985,Sanderetal1988}. Philip \cite{Philip1969} derived an exact solution for infiltration problems expressed in the form $z(\theta,t)$, that is, the depth where the water content takes specified values at given time points $t$. Philip's solution has been used in \cite{Warricketal1985} to construct a table of coefficients which allow the computation of $z(\theta,t)$ for three different $\theta$ and arbitrary $t$. The solution verifies the dimensionless form the $\theta$-based Richards equation
\begin{equation}\label{eqTheta}
\frac{\partial \Theta}{\partial t}=\frac{\partial}{\partial z}\left[D(\Theta)\frac{\partial}{\partial z}(\Theta)\right]-\frac{dK(\Theta)}{d\Theta}\frac{\partial\Theta)}{\partial z},
\end{equation}
where $z$ is positive downward, $D(\Theta)=K(\Theta)d\psi/d\Theta$, and $K(\Theta)$  is given according to the van Genuchten-Mualem model by the upper branch of (\ref{K}). Such analytical solutions have been used in \cite{Phoonetal2007,Caviedesetal2013} to verify various one-dimensional numerical schemes based on finite volume and finite element approaches. In order to test the GRW $L$-scheme (\ref{eq8}-\ref{eq10}), we solve the same infiltration problem (soil column 100 cm deep, constant unsaturated initial water content $\theta_i$, and infiltration imposed by $\psi=0$ on the upper boundary). We use a van Genuchten-Mualem parameter $n=1.5$ together with the parameters of the hypothetical loam soil used in \cite{Warricketal1985}: $K_{sat}=6\cdot 10^{-4}$ cm/s, $\theta_{sat}=0.45$, $\theta_{res}=0.1$, $\theta_{i}=0.17$, $\alpha=0.01$ cm$^{-1}$. The pressure corresponding to the initial water content is obtained by (\ref{theta}), $\psi(\theta_{i})=-24.87$ cm. The computations are carried out with $\Delta z=1$ cm, $\Delta t$ between $9.26\cdot 10^{-4}$ h and $5.23\cdot 10^{-4}$ h, $L=0.2$, and the convergence is achieved after a number of 15 to 160 iterations ($\varepsilon_a=\varepsilon_r=5\cdot 10^{-6}$). The analytical solutions $\tilde{z}(\theta,t)$ for $\theta=$0.24, 0.31, and 0.38 at successive times between 0.5 h and 2 h are obtained with the coefficients for $n=1.5$ and $\Theta(\theta_i)=0.2$ given in \cite[Table 3]{Warricketal1985}. The GRW results $z(\theta,t)$ for the same $\theta$ and $t$ are obtained by linear interpolation of the numerical solution $\theta(z,t)$. Relative errors $(z-\tilde{z})/\tilde{z}$ of the numerical solution $z(\theta,t)$ with respect to the analytical solution $\tilde{z}(\theta,t)$ are shown in Table~\ref{table_Warrick}.

\begin{table}[!ht]
\centering
\caption{Relative errors of the GRW solution $z(\theta,t)$\\ with respect to the exact solution from \cite{Warricketal1985}.}
\label{table_Warrick}
\begin{tabular}{ c c c c c c c c c c}
\hline
\vspace{-0.3cm}\\
   & $t$ (h) & $\theta=0.24$ & $\theta=0.31$ & $\theta=0.38$\\
  \hline
   & 0.5 & 5.31e-02  & 5.31e-02  & 5.69e-02 &\\
   & 1.0 & -2.46e-03 & 2.175e-02 & 5.10e-02 &\\
   & 1.5 & -5.70e-02 & -1.41e-02 & 4.55e-02 &\\
   & 2.0 & -9.69e-02 & -3.88e-02 & 4.91e-02 &\\
  \hline
\end{tabular}
\end{table}

An exact solution for constant flux infiltration with dry initial condition $\Theta(z,0)=0$ has been derived in \cite{Sanderetal1988} and further used to verify the numerical solution provided by a pressure formulation of the Richards equation \cite{Watsonetal1995}. The solution solves Eq.(\ref{eqTheta}) with coefficient given by Fujita's model \cite{Fujita1952},
\[
D(\Theta)=D_0/(1-v\Theta)^2 ,
\]
where $D_0$ and $v$ are positive constants. Since $\Theta(z,0)=0$ implies $\psi(\Theta(z,0))=\infty$ as initial condition for the numerical scheme in pressure formulation, the singularity was avoided in \cite{Watsonetal1995} by considering $\Theta(z,0)=3.4483\cdot 10^{-6}$ as a numerical simulation parameter. As for the GRW scheme (\ref{eq8}-\ref{eq10}), we would have $K(\Theta(z,0))=0$ and, according to (\ref{eq7}), the condition $r^{s}_{i\pm 1/2,k}\le 1/2$ implies $\Delta t=\infty$, for finite $\Delta z$. Using the same initial $\psi$ as in \cite{Watsonetal1995} requires a very fine discretization which would slow down considerably the computation. Therefore, we opt for the direct approach of solving (\ref{eqTheta}) as a diffusion equation with drift coefficient defined by $V(\Theta)=dK(\Theta)/d\Theta$. The latter will be computed analytically from the parameterization $K(\Theta)$ used in \cite{Watsonetal1995}.

Proceeding as in Section~\ref{1Dalg}, we start with a forward-time centered-space finite difference discretization of Eq.~(\ref{eqTheta}),
\begin{align*}
\Theta_{i,k+1}-\Theta_{i,k1}=
&+\frac{\Delta t}{{\Delta z}^2}\left[D(\Theta_{i+1/2,k})(\Theta_{i+1,k}-\Theta_{i,k})-
D(\Theta_{i-1/2,k})(\Theta_{i,k}-\Theta_{i-1,k})\right]\nonumber\\
&-\frac{\Delta t}{{2\Delta z}}V_{i,k}(\Theta_{i+1,k}-\Theta_{i-1,k}),
\end{align*}
we approximate the solution by a distribution of $\mathcal{N}$ particles on a regular lattice, $\Theta_{i,k}\approx n_{i,k}/\mathcal{N}$, and end up with
\begin{align}\label{eqTbgrw1}
n_{i,k+1}=&[1-(r_{i+1/2,k}+r_{i-1/2,k})]n_{i,k}\nonumber\\
&+\frac{1}{2}(r_{i+1/2,k}-v_{i,k})n_{i+1,k}+\frac{1}{2}(r_{i-1/2,k}+v_{i,k})n_{i-1,k}.
\end{align}
The dimensionless parameters in Eq.~(\ref{eqTbgrw1}) are given by
\[
r_{i\pm 1/2,k}=\frac{2\Delta t}{{\Delta z}^2}D_{i\pm1/2,k},\;\; v_{i,k}=\frac{\Delta t}{\Delta z},\;\;
r_{i\pm 1/2,k}\le 1,\;\; |v_{i,k}|\le r_{i\pm 1/2,k}.
\]
Equation~(\ref{eqTbgrw1}) sums up contributions of random walkers jumping on the lattice according to the rule
\begin{equation}\label{eqTbgrw2}
n_{j,k}=\delta n_{j\mid j,k}+\delta n_{j-1\mid j,k}+\delta n_{j+1\mid j,k},
\end{equation}
which defines a biased global random walk algorithm (BGRW) \cite[Sect. 3.3.3]{Suciu2019}. The numbers of particles $\delta n$ in (\ref{eqTbgrw2}) are binomial random variables determined by the same procedure as in Section~\ref{1Dalg} and their ensemble averages verify
\[
\overline{\delta n_{j\mid j,k}}=[1-(r_{i+1/2,k}+r_{i-1/2,k})]\overline{n_{i,k}},\;\; \overline{\delta n_{j\pm 1\mid j,k}}=\frac{1}{2}(r_{i\pm1,k}\mp v_{i,k})\overline{n_{i,k}}.
\]

Following \cite{Watsonetal1995}, we set on the top boundary the constant flux condition $Q=q/(\theta_{sat}-\theta_{res})=0.2759$ cm/min, with $\theta_{sat}=0.35$, $\theta_{res}=0.06$, and consider the constant parameters $D_0=2.75862$ cm$^2$/min and $v=0.85$ of the Fujita's model. The BGRW results for the final time $T=0.3625$ min, obtained with  $\Delta z=10^{-2}$ cm and $\Delta t$ between $1.51\cdot 10^{-5}$ min and $1.23\cdot 10^{-5}$ min, are compared in Table~\ref{table_Sander} with the analytical solution presented in \cite[Table 1]{Watsonetal1995}.

\begin{table}[!ht]
\centering
\caption{GRW solution $\theta(z,t)$ compared to\\ the analytical solution from \cite{Sanderetal1988}.}
\label{table_Sander}
\begin{tabular}{ c c c c c c c c c c}
\hline
\vspace{-0.3cm}&\\
   & $z$ (cm) & $\tilde{\theta}(z,t)$ & $\theta(z,t)$ & $(\theta-\tilde{\theta})/\tilde{\theta}$\\
  \hline
   & 0    & 0.0907  & 0.0929 & 2.53e-02 &\\
   & -0.2 & 0.0861  & 0.0884 & 2.68e-02 &\\
   & -0.4 & 0.0819  & 0.0842 & 2.70e-02 &\\
   & -0.6 & 0.0782  & 0.0802 & 2.60e-02 &\\
   & -0.8 & 0.0748  & 0.0766 & 2.40e-02 &\\
   & -1.0 & 0.0719  & 0.0734 & 2.12e-02 &\\
   & -2.0 & 0.0631  & 0.0635 & 6.40e-03 &\\
  \hline
\end{tabular}
\end{table}

The tests for unsaturated one-dimensional flows presented above are completed in Section~\ref{1Dvalid_flow_transp} by convergence investigations and estimations of convergence order of the GRW algorithms for fully coupled nonlinear flow and transport problems for saturated/unsaturated porous systems.

\section{Two-dimensional GRW solutions}
\label{2D}

\subsection{Two-dimensional GRW algorithm for flow in soils and aquifers}
\label{2DalgF}

In two spatial dimensions the pressure head $\psi(x,z,t)$ satisfies the equation
\begin{equation}\label{eq11}
\frac{\partial}{\partial t}\theta({\psi})-\nabla\cdot\left[K(\theta(\psi)\nabla(\psi+z)\right]=0.
\end{equation}

The two-dimensional GRW algorithm on regular staggered grids ($\Delta x=\Delta z$) which approximates the solution of (\ref{eq11}) by computational particles, $\psi\approx n a /\mathcal{N}$, is constructed similarly to (\ref{eq8}-\ref{eq10}). The solution at iteration $s+1$ is obtained by gathering particles from neighboring sites according to
\begin{eqnarray}\label{eq12}
n_{i,j,k}^{s+1}
&=&\left[1-\left(r_{i+1/2,j,k}^{s}+r_{i-1/2,j,k}^{s}+r_{i,j+1/2,k}^{s}+r_{i,j-1/2,k}^{s}\right)\right]n_{i,j,k}^{s}\nonumber\\
& & + r_{i+1/2,j,k}^{s}n_{i+1,j,k}+r_{i-1/2,j,k}^{s}n_{i-1,j,k}\nonumber\\
& & + r_{i,j+1/2,k}^{s}n_{i,j+1,k}+r_{i,j-1/2,k}^{s}n_{i,j-1,k}+\left\lfloor \mathcal{N}f^{s}\right\rfloor,
\end{eqnarray}
where the source term is defined as $f^{s}=\left(r_{i,j+1/2,k}^{s}-r_{i,j-1/2,k}^{s}\right)\Delta z
-\left[\theta(n_{i,j,k}^{s})-\theta(n_{i,j,k-1})\right]/L$.
The two-dimensional GRW rule which at time $k$ moves particles from sites $(l,m)$ to neighboring sites $(l\mp 1,m\mp 1)$ reads as follows,
\begin{equation}\label{eq13}
n_{l,m,k}^{s} = \delta n_{l,m|l,m,k}^{s}
+\delta n_{l-1,m|l,m,k}^{s}+\delta n_{l+1,m|l,m,k}^{s}
+\delta n_{l,m-1|l,m,k}^{s}+\delta n_{l,m+1|l,m,k}^{s}.
\end{equation}
For consistency with (\ref{eq12}), the numbers of particles $\delta n^{s}$ verify in the mean
\begin{eqnarray}\label{eq14}
&&\overline{\delta n_{l,m|l,m,k}^{s}}
=\left[1-\left(r_{l-1/2,m,k}^{s}+r_{l+1/2,m,k}^{s}+r_{l,m-1/2,k}^{s}+r_{l,m+1/2,k}^{s}\right)\right]\overline{n_{l,m,k}^{s}}\nonumber\\
&&\overline{\delta n_{l\mp 1,m|l,m,k}^{s}}=r_{l\mp 1/2,m,k}^{s}\overline{n_{l,m,k}^{s}}\nonumber\\
&&\overline{\delta n_{l,m\mp 1|l,m,k}^{s}}=r_{l,m\mp 1/2,k}^{s}\overline{n_{l,m,k}^{s}}.
\end{eqnarray}
The parameters $r_{l\mp 1/2,m,k}^{s}$ and $r_{l,m\mp 1/2,k}^{s}$, defined by relations similar to (\ref{eq7}), are dimensionless positive real numbers. They represent biased jump probabilities on the four allowed spatial directions of the GRW lattice and are constraint by the first relation (\ref{eq14}) such that their sum be less or equal to one. A sufficient condition would be that each of them verifies $r\leq 1/4$.

The binomial random variables variables $\delta n$ are approximated in the same way as in the one-dimensional case. By giving up the particle indivisibility, one obtains deterministic GRW algorithms which represent the solution $n$ by real numbers and use the unaveraged relations (\ref{eq14}) for the computation of the $\delta n$ terms. In the following we use this deterministic implementation of the GRW algorithm to compute flow solutions for unsaturated/saturated porous media.

\begin{remark}\label{remFsat2D}
After disregarding the index $k$ and letting $L=1/a$, $\theta=const$, the algorithm (\ref{eq12}-\ref{eq14}) becomes a transient scheme to solve the equation governing flows in saturated porous media \cite{Alecsaetal2020,Suciu2020} (see also Remark~\ref{remFsat}).
\end{remark}

\subsection{GRW algorithms for two-dimensional fully coupled flow and surfactant transport}
\label{2DalgFT}
Let the pressure $\psi(x,z,t)$ and the concentration $c(x,z,t)$ solve the equations of the following model of fully coupled flow and surfactant transport in unsaturated/saturated porous media \cite{Knabneretal2003,Illiano2020},
\begin{equation}\label{eqFT01}
\frac{\partial}{\partial t}\theta(\psi,c)-\nabla\cdot\left[K(\theta(\psi,c)\nabla(\psi+z)\right]=0,
\end{equation}

\begin{equation}\label{eqFT02}
\frac{\partial}{\partial t}\left[\theta(\psi,c)c\right]-\nabla\cdot\left[D\nabla c -\mathbf{q}c\right]=R(c),
\end{equation}
where $\mathbf{q}=-K(\theta(\psi,c)\nabla(\psi+z)$ is the water flux (Darcy velocity) and $R(c)$ is a nonlinear reaction term. Equations~(\ref{eqFT01}-\ref{eqFT02}) are coupled in both directions through the nonlinear functions $\theta(\psi,c)$ and $\theta(\psi,c)c$. The pressure equation (\ref{eqFT01}) is solved with the GRW $L$-scheme described in the previous subsection, with a slight modification due to the dependence of $\theta$ on both $\psi$ and $c$.
New algorithms are needed instead to solve the coupled, nonlinear transport equation (\ref{eqFT02}).

\subsubsection{Biased GRW algorithm for transport problems}
\label{2DalgFT_BGRW}

To derive a GRW algorithm for the transport equation, we start with a backward-time central-space finite difference scheme for Eq.~(\ref{eqFT02}). Considering a diagonal diffusion tensor with constant 
components $D_1$ and $D_2$, and denoting by $U$ and $V$ the components of the Darcy velocity along the horizontal axis $x$ and the vertical axis $z$, by $\Delta t$ the time step, and by $\Delta x$ and $\Delta z$ the spatial steps, the scheme reads as
\begin{eqnarray}\label{eqFT1}
&&\theta(\psi_{i,j,k},c_{i,j,k})c_{i,j,k}-\theta(\psi_{i,j,k-1},c_{i,j,k-1})c_{i,j,k-1}=\nonumber\\ &&-\frac{\Delta t}{2\Delta x}\left(U_{i+1,j,k}c_{i+1,j,k}-U_{i-1,j,k}c_{i-1,j,k}\right)-\frac{\Delta t}{2\Delta z}\left(V_{i,j+1,k}c_{i,j+1,k}-V_{i,j-1,k}c_{i,j-1,k}\right)\nonumber\\
&&+\frac{D_1\Delta t}{\Delta x ^2}\left(c_{i+1,j,k}-2c_{i,j,k}+c_{i-1,j,k}\right)+\frac{D_2\Delta t}{\Delta z ^2}\left(c_{i,j+1,k}-2c_{i,j,k}+c_{i,j-1,k}\right)+R(c_{i,j,k})\Delta t=\nonumber\\
&&-\left(\frac{2D_1\Delta t}{\Delta x ^2}+\frac{2D_2\Delta t}{\Delta z ^2}\right)c_{i,j,k}\nonumber\\
&&+\left(\frac{D_1\Delta t}{\Delta x ^2}-\frac{\Delta t}{2\Delta x}U_{i+1,j,k}\right)c_{i+1,j,k}+\left(\frac{D_1\Delta t}{\Delta x ^2}+\frac{\Delta t}{2\Delta x}U_{i-1,j,k}\right)c_{i-1,j,k}\nonumber\\
&&+\left(\frac{D_2\Delta t}{\Delta z ^2}-\frac{\Delta t}{2\Delta z}V_{i,j+1,k}\right)c_{i,j+1,k}+\left(\frac{D_2\Delta t}{\Delta z ^2}+\frac{\Delta t}{2\Delta z}V_{i,j-1,k}\right)c_{i,j-1,k}+R(c_{i,j,k})\Delta t.
\end{eqnarray}
Next, similarly to the scheme for the flow equation, we add a regularization term $L(c^{s+1}_{i,j,k}-c^{s}_{i,j,k})$ in Eq.~(\ref{eqFT1}), define the dimensional parameters
\begin{equation}\label{eqFT2}
r_{x}=\frac{2D_1\Delta t}{L\Delta x ^2},\;\; r_{z}=\frac{2D_2\Delta t}{L\Delta z ^2},\;\; u_{i\pm 1,j,k}^{s}=\frac{\Delta t}{L\Delta x}U_{i\pm 1,j,k}^{s},\;\; v_{i,j\pm 1,k}^{s}=\frac{\Delta t}{L\Delta z}V_{i,j\pm 1,k}^{s},
\end{equation}
approximate the concentration by the density of the number of computational particles, $c^{s}_{i,j,k}\approx n^{s}_{i,j,k}/\mathcal{N}$, and finally we obtain
\begin{align}\label{eqFT3}
n_{i,j,k}^{s+1}
=&\left[1-\left(r_{x}+r_{z}\right)\right]n_{i,j,k}^{s}\nonumber\\
& + \frac{1}{2}\left(r_{x}-u_{i+ 1,j,k}^{s}\right)n_{i+1,j,k}^{s}+ \frac{1}{2}\left(r_{x}+u_{i-1,j,k}^{s}\right)n_{i-1,j,k}^{s}\nonumber\\
& + \frac{1}{2}\left(r_{z}-v_{i,j+ 1,k}^{s}\right)n_{i,j+1,k}^{s}+ \frac{1}{2}\left(r_{z}+v_{i,j-1,k}^{s}\right)n_{i,j-1,k}^{s} +\left\lfloor \mathcal{N}g^{s}\right\rfloor,
\end{align}
where $g^{s}=R(n_{i,j,k}^{s})\Delta t/L-\left[\theta(\psi_{i,j,k}^{s},n_{i,j,k}^{s})n_{i,j,k}^{s}
-\theta(\psi_{i,j,k-1},n_{i,j,k-1})n_{i,j,k-1}\right]/L$, with $\psi$ approximated by the distribution of particles in the flow solver for Eq.~(\ref{eqFT01}). Note that the definition of the dimensionless numbers (\ref{eqFT2}) implies that the parameter $L$ has to be a dimensionless number as well.

The contributions to $n_{i,j,k}^{s+1}$ in Eq.~(\ref{eqFT3}) are obtained with the BGRW algorithm
\begin{equation}\label{eqFT4}
n_{l,m,k}^s=\delta n_{l,m\mid l,m,k}^s + \delta n_{l-1,m\mid l,m,k}^s+\delta n_{l+1,m\mid l,m,k}^s + \delta n_{l,m-1\mid l,m,k}^s+\delta n_{l,m+1\mid l,m,k}^s,
\end{equation}
where, for consistency with the finite difference scheme~(\ref{eqFT3}), the quantities $\delta n$ verify in the mean
\begin{align}\label{eqFT5}
\overline{\delta n_{l,m\mid l,m,k}^s}=\left[1-\left(r_{x}+r_{z}\right)\right]\mbox{
}\overline{n_{i,j,k}^s},\;\;\;
&\overline{\delta n_{l\pm 1,m\mid l,m,k}^s}=\frac{1}{2}(r_{x}\mp u_{l,m,k}^{s})\overline{n_{l,m,k}^s},\nonumber\\
&\overline{\delta n_{l,m\pm 1\mid l,m,k}^s}=\frac{1}{2}(r_{z}\mp
v_{l,m,k}^{s})\overline{n_{l,m,k}^s}.
\end{align}
The binomial random variables $\delta n$ used in the BGRW algorithm are approximated similarly to the algorithms described in the previous sections, by summing up to unity reminders of multiplication and floor operations. A deterministic BGRW algorithm can be obtained, similarly to the flow solver presented in Section~\ref{2DalgF} above, by giving up the particle's indivisibility and using the un-averaged relations (\ref{eqFT5}). However, for the computations presented in the next section, we use a randomized implementation of the BGRW algorithm.

As follows from (\ref{eqFT5}), the BGRW algorithm is subject to the
following restrictions
\begin{equation}
r_{x}+r_{z}\leq1,\mbox{  }\left\vert u_{l,m,k}^{s}\right\vert \leq
r_{x},\mbox{  }\left\vert v_{l,m,k}^{s}\right\vert \leq r_{z}.\label{eqFT6}
\end{equation}

\begin{remark}\label{remPe} The constraints (\ref{eqFT6}) impose a limitation on the maximum allowable value of the local P\'{e}clet number. Assume a constant flow velocity $-V$ and a constant diffusion coefficient $D$. Then, according to (\ref{eqFT6}) and (\ref{eqFT2}), the condition $v\le r$ implies $\text{P\'{e}}=V\Delta z/D\le 2$.
\end{remark}

\begin{remark}\label{remTdec}
Taking into account that the iterations start with $n_{i,j,k-1}$, setting $L=1$, $\theta=1$, and dropping the superscripts $s$, the relation (\ref{eqFT3}) becomes
\begin{align}\label{eqFT3a}
n_{i,j,k}
=&\left[1-(r_x+r_z)\right]n_{i,j,k-1}\nonumber\\
& + \frac{1}{2}\left(r_{x}-u_{i+ 1,j,k-1}\right)n_{i+1,j,k-1}+ \frac{1}{2}\left(r_{x}+u_{i-1,j,k-1}\right)n_{i-1,j,k-1}\nonumber\\
& + \frac{1}{2}\left(r_{z}-v_{i,j+ 1,k-1}\right)n_{i,j+1,k-1}+ \frac{1}{2}\left(r_{z}+v_{i,j-1,k-1}\right)n_{i,j-1,k-1} +\left\lfloor \mathcal{N}R(n_{i,j,k-1})\Delta t\right\rfloor.
\end{align}
Relation (\ref{eqFT3a}), together with (\ref{eqFT4}-\ref{eqFT6}), define a BGRW algorithm for (decoupled) reactive transport described by Eq.~(\ref{eqFT02}) with $\theta(\psi,c)=1$.
\end{remark}

\subsubsection{Unbiased GRW algorithm for transport problems}
\label{2DalgFT_GRW}

The unbiased GRW algorithm is obtained by globally moving groups of particles according to the rule
\begin{align}
n_{i,j,k}^{s}=&\;\delta n_{i+u_{i,j,k}^{s},j+v_{i,j,k}^{s}\mid i,j,k}^s\label{eqFT7}\\
&+ \delta n_{i+u_{i,j,k}^{s}+d,j+v_{i,j,k}^{s}\mid i,j,k}^s+\delta n_{i+u_{i,j,k}^{s}-d,j+v_{i,j,k}^{s}\mid i,j,k}^s\nonumber \\
&+ \delta n_{i+u_{i,j,k}^{s},j+v_{i,j,k}^{s}+d\mid i,j,k}^s+\delta n_{i+u_{i,j,k}^{s},j+v_{i,j,k}^{s}-d\mid i,j,k}^s,\nonumber
\end{align}
where $d$ is a constant amplitude of diffusion jumps and the dimensionless variables $r_x$, $r_z$, $u$ and $v$ are defined similarly to (\ref{eqFT2}) by
\begin{equation}\label{eqFT2a}
r_x=\frac{2D_1\Delta t}{L(d\Delta x )^2},\;\; r_z=\frac{2D_2\Delta t}{L(d\Delta z )^2},\;\; u_{i,j,k}^{s}=\left\lfloor\frac{\Delta t}{L\Delta x}U_{i,j,k}^{s}+0.5\right\rfloor,\;\; v_{i,j,k}^{s}=\left\lfloor\frac{\Delta t}{L\Delta z}V_{i,j,k}^{s}+0.5\right\rfloor.
\end{equation}

The particles distribution is updated at every time step by
\begin{equation}
n_{l,m,k}^{s+1}=\delta n_{l,m,k}^{s}+\sum_{i\neq l,j\neq m}\delta n_{l,m\mid i,j,k}^{s} +\left\lfloor \mathcal{N}g^{s}\right\rfloor.\label{eqFT8}
\end{equation}

The averages over GRW runs of the terms from (\ref{eqFT7}) are now
related by
\begin{align}
&\overline{\delta n_{i+u_{i,j,k}^{s},j+v_{i,j,k}^{s}\mid i,j,k}^{s}}=\left[1-\left(r_{x}+r_{z}\right)\right]\mbox{
}\overline{n_{i,j,k}^{s}},\nonumber\\
&\overline{\delta n_{i+u_{i,j,k}^{s}\pm d,j+v_{i,j,k}^{s}\mid i,j,k}^{s}}=\frac{r_{x}}{2}\hspace{0.1cm}\overline{n_{i,j,k}^{s}},\nonumber\\
&\overline{\delta n_{i+u_{i,j,k}^{s},j+v_{i,j,k}^{s}\pm d\mid i,j,k}^{s}}=\frac{r_{z}}{2}\hspace{0.1cm}\overline{n_{i,j,k}^{s}}.\label{eqFT9}
\end{align}
Comparing with the BGRW relations (\ref{eqFT5}), we remark that (\ref{eqFT2a}) defines unbiased jump probabilities $r_x/2$ and $r_y/2$ on the two spatial directions.

The unbiased GRW algorithm for decoupled transport is obtained by letting $L=1$ and dropping the superscripts $s$ (see also Remark~\ref{remTdec}).

The binomial random variables $\delta n$ used in the unbiased GRW algorithm are approximated by the procedure used for the flow solver and for the BGRW algorithm presented in the previous subsection. For fixed space steps, the time step is chosen such that the dimensionless parameters $u_{i,j,k}^{s}$ and $v_{i,j,k}^{s}$ take integer values larger than unity which ensure the desired resolution of the velocity components \cite[Sect. 3.3.2.1]{Suciu2019}. Further, the jumps' amplitude $d$ is chosen such that the jump probabilities verify the constraint $r_{x}+r_{z}\leq 1$, imposed by the first relation (\ref{eqFT9}).

The unbiased GRW, as well as the BGRW algorithm introduced in Section~\ref{2DalgFT_BGRW} above, have been tailored to solve problems with constant diffusion coefficients, as those considered in Sections~\ref{2DvalidFT}~and~\ref{2D_FlowTransp_aquifers} below. In case of diagonal diffusion tensors with space-time variable coefficients $D_1$ and $D_2$, the algorithms for the transport problem are straightforwardly obtained by assigning to $r_x$ and $r_z$ superscripts $s$ and appropriate subscripts $i,j,k$.

\section{Validation of the two-dimensional GRW algorithms}
\label{2Dvalid}

\subsection{GRW flow solutions}
\label{2DvalidF}

For the beginning, we conduct verification tests of the GRW flow code by comparisons with an analytical solution and compute numerical estimates of the order of convergence. The results are further compared with those obtained by a TPFA code implemented in the MRST software \cite{mrst,Illiano2020}. The two codes are tested by solving a problem with manufactured solution previously considered in \cite{Radu2014}. The domain is the unit square $[0,1]\times [0,1]$ and the final time is  $T=1$. The manufactured solution for the pressure head $\psi_m$ is given by
\begin{equation}\label{eqFsol}
\psi_m(x,z,t) = - t\ x\ (x-1)\ z\ (z-1)\ -\ 1.
\end{equation}
The water content $\theta$ and the conductivity $K$ are expressed as
\begin{equation}\label{eqFparam}
\theta(\psi) = \frac{1}{1-\psi}\;, \quad K(\theta(\psi)) = \psi^{2}\;.
\end{equation}
The analytical solution (\ref{eqFsol}) defines the boundary and initial conditions and induces a source term $f$, computed analytically from Eq.~(\ref{eq11}) with parameters given by the expressions (\ref{eqFparam}).

We start the computations on a uniform mesh with $\Delta x = \Delta z = 0.1$ and halve the mesh size step three times successively. The accuracy of the numerical solutions, at the final time $t=T$, is quantified by the $L^2$ norm $\varepsilon_l=\|\psi^{(l)}-\psi_m\|$, $l=1,\ldots,4$, where $l=1$ corresponds to the original mesh. The estimated order of convergence (EOC) that describes the decrease of the error in logarithmic scale is computed according to
\begin{equation}\label{eq16}
EOC=\log\left(\frac{\varepsilon_l}{\varepsilon_{l+1}}\right)/\log(2), \quad l=1,\ldots,3.
\end{equation}

The computations with the TPFA code start with a time step $\Delta t = 0.1$ which is also halved at each refinement of the mesh. The parameters of the convergence indicator (\ref{eq5}) are set to $\varepsilon_a=10^{-6}$ and $\varepsilon_r=0$. Finally, the linearization parameter $L$ is set equal to $1/2$ and the convergence of the $L$-scheme is achieved after circa 100 iterations per time step, independently of the mesh size.

In the GRW computations we use the same spatial refinement of the grid and tolerances $\varepsilon_a$ and $\varepsilon_r$ as above but, according to (\ref{eq7}), we have to use adaptive time steps $\Delta t = \mathcal{O}(\Delta z^{1/2})$ (see discussion in Section~\ref{1Dvalid}). The convergence criterion (\ref{eq5}) is already fulfilled by the GRW $L$-scheme with parameter $L=1$ for numbers of iterations increasing from $s=2$ to $s=5$ as the space step decreases. The accuracy $\varepsilon_l$ instead is strongly influenced by $L$. For $L<800$ the $\varepsilon_l$ values may increase with the refinement of the mesh, leading to negative EOC, that is, the GRW solution does not converge to the exact solution $\psi_m$. However, it is found that the increase of $\varepsilon_l$ is prevented by using a sufficiently large parameter $L$.

The results presented in Table \ref{tab:Richards_2D_Test_MRST} indicate the convergence of order 1 in space for TPFA and of order 2 for the GRW solutions. The higher order of convergence also leads to much smaller errors of the GRW code after the first refinement of the mesh.

\begin{table}[h!]
\begin{center}
\caption{Estimated order of convergence of the TPFA and GRW flow solvers.}
\label{tab:Richards_2D_Test_MRST}
\begin{tabular}{c c c c c c c c}
\hline
&$\varepsilon_1$ & EOC & $\varepsilon_2$ & EOC & $\varepsilon_3$ & EOC & $\varepsilon_4$\\
\hline
TPFA        & 8.45e-03 & 0.94 & 4.40e-03 & 0.97 & 2.25e-03 & 0.97 & 1.15e-03\\
\hline
GRW (L=800) & 7.20e-03 & 2.24 & 1.52e-03 & 3.21 & 1.65e-04 & 0.50 & 1.17e-04 \\ \hline
GRW (L=1000)& 9.24e-03 & 2.22 & 1.99e-03 & 2.83 & 2.80e-04 & 1.66 & 8.84e-05 \\ \hline
GRW (L=1200)& 8.89e-03 & 2.23 & 1.90e-03 & 2.80 & 2.72e-04 & 2.14 & 6.16e-05 \\
\hline
\end{tabular}
\end{center}
\end{table}

Further, we solve the benchmark problem from \cite[Sect. 4.2]{ListandRadu2016}, which describes the recharge of a groundwater reservoir from a drainage trench in a two-dimensional geometry. The groundwater table is fixed by a Dirichlet boundary condition on the right hand side. The drainage process is driven by a Dirichlet boundary condition changing in time on the upper boundary of $\Omega$.

The precise structure of the domain is defined by
\begin{equation*}
\begin{split}
\Omega & = (0,2)\times (0,3),\\
\Gamma_{D_1} &= \{(x,z) \in \partial \Omega\ |\ x \in [0,1] \wedge z = 3\},\\
\Gamma_{D_2} &= \{(x,z) \in \partial \Omega\ |\ x =2 \wedge z \in [0,1]\},\\
\Gamma_D &=  \Gamma_{D_1} \cup \Gamma_{D_2},\\
\Gamma_N &= \partial \Omega \setminus \Gamma_D.
\end{split}
\end{equation*}
The Dirichlet and Neumann boundary conditions on $\Gamma_D$ and $\Gamma_N$, respectively, as well as the initial condition consisting of hydrostatic equilibrium are specified as follows:
\begin{equation*}
\begin{split}
&\psi(x,z,t) = \begin{cases}
-2 + 2.2t/\Delta t_D, &\quad\ \text{on}\ \Gamma_{D_1}, T\leq \Delta t_D,\\
0.2, &\quad\ \text{on}\ \Gamma_{D_1}, T> \Delta t_D,\\
1-z, &\quad\ \text{on} \ \Gamma_{D_2},
\end{cases}\\
&-K(\theta(\psi(x,z,t))\nabla(\psi(x,z,t)+z)\cdot \mathbf{n} \ = \ 0 , \quad\text{on}\ \Gamma_N,\\
&\psi(x,z,0)\ = \ 1-z,  \quad\text{on}\ \Omega,\\
\end{split}
\end{equation*}
where $\mathbf{n}$ represents the outward pointing normal vector.

We consider here two sets of soil parameters, presented in Table \ref{table:parametersRichards}, which correspond to a silt loam and a Beit Netofa clay, respectively.
\begin{table}[!ht] \centering
 \caption{Simulation parameters.}
 \label{table:parametersRichards}
 \begin{tabular}{ c c c }
  \hline
  & Silt loam & Beit Netofa clay \\
  \hline
Vam Genuchten parameters:&& \\
$\theta_{sat}$ & 0.396 & 0.446\\
$\theta_{res}$ & 0.131 & 0\\
$\alpha$ & 0.423 & 0.152\\
$n$ & 2.06 & 1.17\\
$K_{sat}$ & $4.96 \cdot 10^{-2}$ & $8.2\cdot 10^{-4}$\\
\hline
Time parameters: & &\\
$\Delta t_D$ & 1/16 & 1\\
$\Delta t$ & 1/48 & 1/3\\
$T$ & 3/16 & 3\\
\hline
 \end{tabular}
\end{table}

The time unit is 1 day and spatial dimensions are given in meters. Furthermore, we consider a regular mesh consisting of 651 nodes (i.e., $\Delta x = \Delta z = 0.1$).

\begin{figure}
\begin{minipage}[t]{0.45\linewidth}\centering
\includegraphics[width=\linewidth]{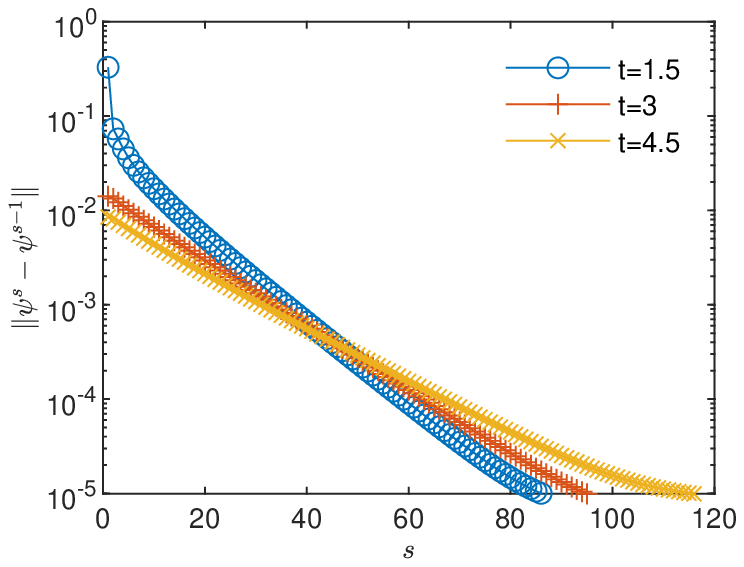}
\caption{\label{fig14}Convergence of the $L$-scheme implementation of the GRW flow solver for the loam soil problem at three time levels (in hours).}
\end{minipage}
\hspace*{0.1in}
\begin{minipage}[t]{0.45\linewidth}\centering
\includegraphics[width=\linewidth]{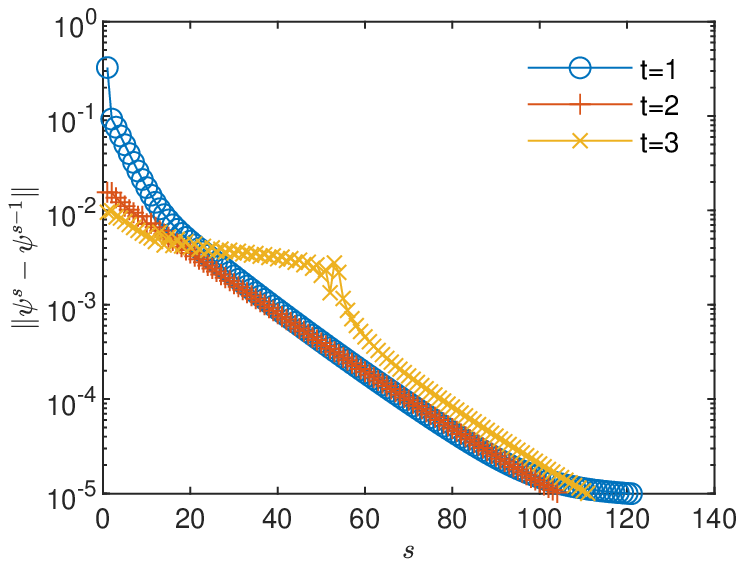}
\caption{\label{fig16}Convergence of the $L$-scheme implementation of the GRW flow solver for the clay soil problem at three time levels (in days).}
\end{minipage}
\end{figure}

By setting the stabilization parameters to $L=0.5$ for loam and for $L=0.12$ for clay, the convergence criterion (\ref{eq5}) with $\varepsilon_a=\varepsilon_r=5\cdot10^{-6}$ is fulfilled after about 120 iterations of the GRW $L$-scheme, for both soil models (Figs.~\ref{fig14}~and~\ref{fig16}).
The computational orders of convergence of the $L$-scheme indicate linear convergence for both loam and clay soil models (see Appendix~\ref{appA}).
The results shown in Figs.~\ref{fig13}~and~\ref{fig15} are as expected for this benchmark problem (see \cite{Schneid2000,ListandRadu2016}): the drainage process in the clay soil is much slower, so that the pressure distribution after three days is similar to that established in the loam soil after 4.5 hours.

The results obtained with the TPFA $L$-scheme, with $L=1$ for both soil models, are used as reference to compute the relative errors $\varepsilon_\psi$, $\varepsilon_\theta$, $\varepsilon_{q_x}$, and $\varepsilon_{q_z}$ shown in Table \ref{table_relative_errors}. One remarks that $\varepsilon_\psi$ and $\varepsilon_\theta$ are close to the corresponding errors for the one-dimensional case presented in Table~\ref{table1}, but $\varepsilon_{q_x}$ and $\varepsilon_{q_z}$ are one order of magnitude larger
than $\varepsilon_{q}$ in shown in Table~\ref{table1}. A possible explanation could be the occurrence of the numerical diffusion in the flow TPFA code (see discussion at the end of Section~\ref{coupledFT} below). The computational times of the GRW code are 1 second and 1.6 seconds for loam and clay cases, respectively. The times of the TPFA runs, on the same computer, are one order of magnitude larger, i.e., 25 seconds and 38 seconds, respectively.

\begin{figure}
\begin{minipage}[t]{0.45\linewidth}\centering
\includegraphics[width=\linewidth]{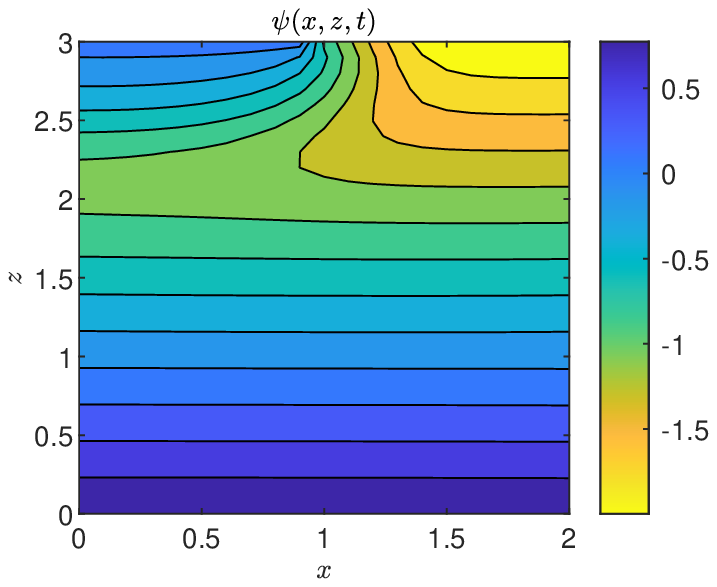}
\caption{\label{fig13}Pressure head solution at $t=4.5$ hours obtained by the GRW code for the benchmark problem of recharge from a drainage trench through a silt loam soil.}
\end{minipage}
\hspace*{0.1in}
\begin{minipage}[t]{0.45\linewidth}\centering
\includegraphics[width=\linewidth]{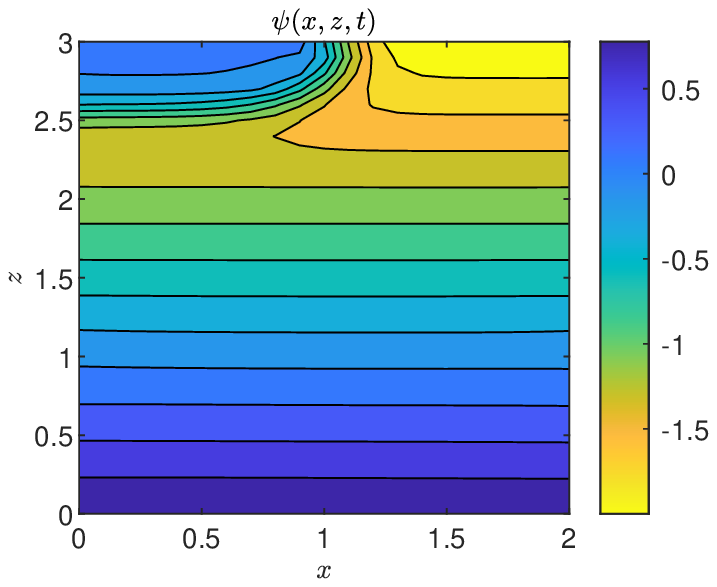}
\caption{\label{fig15}Pressure head solution at $t=3$ days obtained by the GRW code for the benchmark problem of recharge from a drainage trench through a Beit Netofa clay soil.}
\end{minipage}
\end{figure}

\begin{table}[!ht]
\centering
\caption{Comparison of GRW and TPFA\\ solutions of the flow benchmark problem.}
\label{table_relative_errors}
\begin{tabular}{ c c c c c  c c c}
\hline
& $\varepsilon_\psi$  & $\varepsilon_\theta$  & $\varepsilon_{q_x}$ & $\varepsilon_{q_z}$\\ \hline
loam   & 5.73e-02 & 4.00e-03 & 2.30e-01  & 1.04e-01\\
clay   & 5.48e-02 & 6.71e-04 & 4.73e-01  & 1.14e-01\\ \hline
\end{tabular}
\end{table}

\subsection{GRW/BGRW solutions for fully coupled flow and transport problems}
\label{2DvalidFT}

\subsubsection{Code verification tests}
\label{code_verif_FT}

The code verification tests for coupled flow and transport problems are conducted similarly to those for the flow solver presented in the previous subsection, by considering, along with the exact flow solution (\ref{eqFsol}), the exact solution for the concentration field given by
\begin{equation}\label{eqTsol}
c_m(x,z,t) = t\ x\ (x-1)\ z\ (z-1)\ +\ 1.
\end{equation}
After setting $R=0$ and $D=1$, the coupled system of equations (\ref{eqFT01}-\ref{eqFT02}) is solved in the unit square for a total time $T=1$, with source terms, initial conditions, and boundary conditions resulted from  the exact solutions (\ref{eqFsol}) and (\ref{eqTsol}) with a new parameterization given by
\begin{equation}\label{eqFTparam}
\theta(\psi,c) = \frac{1}{1-\psi - c/10}\;, \quad K(\theta(\psi)) = \psi^{2}\;.
\end{equation}

The GRW flow-algorithm (\ref{eq12}-\ref{eq14}), with $\theta$ and $K$ given by (\ref{eqFTparam}), is coupled with the BGRW transport-algorithm (\ref{eqFT3}-\ref{eqFT6}) initialized with $\mathcal{N}=10^{24}$ particles into an alternating splitting scheme \cite{Illiano2020}. The approach alternates iterations of flow and transport solvers until the convergence criterion (\ref{eq5}) with $\varepsilon_a=10^{-6}$ and $\varepsilon_r=0$ is fulfilled by the numerical solutions for both $\psi $ and $c$. In order to highlight the approach to the convergence order 2, the stabilization parameters of the flow and the transport solvers are set to $L_p=L_c=100$. The GRW results presented in Tables~\ref{tab:mrst-grwFTp}~and~\ref{tab:mrst-grwFTc} are compared with results obtained with a TPFA solver applying the same alternating linearized splitting procedure with parameters $L_p = L_c = 1$ which ensure the convergence of order 1.

The GRW flow solver approximates the Darcy velocity by centered differences only in the interior $\Omega$ of the computational domain. Therefore, the velocity $\mathbf{q}|_{\partial\Omega}$, needed to compute the number of biased jumps from the boundary $\partial\Omega$ in the BGRW relation (\ref{eqFT3}) has to be provided in some way. The straightforward approach is to compute the velocity by using an approximate forward finite difference discretization of Darcy's law. Another option is to extend on the boundary the velocity from the first neighboring interior site. Thanks to the manufactured solution (\ref{eqFsol}) on which the code verification test is based, we also have the exact velocity computed analytically. The latter allows accuracy assessments for the above approximations. We note that the GRW results for the pressure solver obtained with analytical, approximate, and extend $\mathbf{q}|_{\partial\Omega}$ are identical in the precision of three significant digits (Table~\ref{tab:mrst-grwFTp}). For the concentration solutions (Table~\ref{tab:mrst-grwFTc}), we note the remarkably good performance of approximate and extended $\mathbf{q}|_{\partial\Omega}$.

\begin{table}[h!]
\begin{center}
\caption{Estimated order of convergence of the TPFA and GRW solvers: pressure solutions.}
\label{tab:mrst-grwFTp}
\begin{tabular}{c c c c c c c c}
\hline
&$\varepsilon_1$ & EOC & $\varepsilon_2$ & EOC & $\varepsilon_3$ & EOC & $\varepsilon_4$\\
\hline
TPFA  & 8.14e-03 & 0.93 & 4.27e-03 & 0.95 & 2.20e-03 & 0.97 &    1.12e-03  \\
\hline
GRW  & 3.71e-03 & 2.02 & 9.18e-04 & 1.94 & 2.40e-04 & 1.45 & 8.78e-05 \\
\hline
\end{tabular}
\end{center}
\end{table}

\begin{table}[h!]
\begin{center}
\caption{Estimated order of convergence of the TPFA and GRW solvers: concentration solutions.}
\label{tab:mrst-grwFTc}
\begin{tabular}{c c c c c c c c}
\hline
&$\varepsilon_1$ & EOC & $\varepsilon_2$ & EOC & $\varepsilon_3$ & EOC & $\varepsilon_4$\\
\hline
TPFA  & 6.26e-03 & 0.83 & 3.52e-03 & 0.89 & 1.90e-03 & 0.91 & 1.01e-03        \\
\hline
GRW (analytical $\mathbf{q}|_{\partial\Omega}$) & 3.92e-03 & 2.00 & 9.78e-04 & 1.83 & 2.74e-04 & 1.05 & 1.32e-04 \\
\hline
GRW (approximate $\mathbf{q}|_{\partial\Omega}$) & 4.72e-03 & 1.99 & 1.19e-03 & 1.85 & 3.29e-04 & 1.17 & 1.46e-04 \\ \hline
GRW ($\mathbf{q}|_{\partial\Omega}$ from $\mbox{int}(\Omega)$) & 5.26e-03 & 2.00 & 1.31e-03 & 1.87 & 3.59e-04 & 1.23 & 1.53e-04 \\
\hline
\end{tabular}
\end{center}
\end{table}

\subsubsection{Estimates of numerical diffusion}
\label{NumDiff}

The small errors shown in Table~\ref{tab:mrst-grwFTc} indicate that the numerical diffusion in solving the transport step of the coupled problem does not play a significant role. This is somewhat expected for the small P\'{e}clet numbers of order P\'{e}=$10^{-2}$ encountered in these computations. But for the numerical setup of the benchmark problem presented in Section~\ref{2DvalidF} and realistic transport parameters P\'{e} can be significantly larger than unity. Therefore we proceed to estimate the numerical diffusion of the codes compared here by following the procedure used in \cite{Raduetal2011}.

We consider the analytical Gaussian solution $c(x,z,t)$ of Eq.~(\ref{eqFT02}) with $\theta=1$, $R=0$, and constant coefficients $D=0.001$ and $V=-0.0331$, corresponding to the Cauchy problem with a Dirac initial concentration pulse located at the coordinates (1,2.1). The constant velocity $V$, oriented downwards along the $z$-axis, is the steady-state solution of the benchmark flow problem from Section~\ref{2DvalidF} with $K=K_{sat}$ corresponding to the loam soil, initial condition $\psi(x,z,0) = 1-z/3$, Dirichlet boundary conditions $\psi(x,0,t)=1$, $\psi(x,3,t)=0$, and no-flow Neumann conditions on the vertical boundaries. The initial condition $c(x,z,0)$ is the same Gaussian function evaluated at $t=1$ and the final time is $T=3$. For decreasing mesh sizes $\Delta x$ and $\text{P\'{e}} =V\Delta x/D$, the number of time steps was restricted by the requirement that the support of the numerical solution does not extend beyond the boundaries $\partial\Omega$ (to mimic diffusion in unbounded domains). The effective diffusion coefficients $D_x$ and $D_z$ are computed from the spatial moments along the $x$- and $z$-directions of the numerical solution (see \cite[Eqs. (38-41)]{Raduetal2011}). The numerical diffusion is estimated by relative errors $\varepsilon_{D_{x}}=|D_x-D|/D$ and $\varepsilon_{D_{z}}=|D_z-D|/D$ averaged over the time interval $[0,T]$. Table~\ref{table_relative_errorsNumDiff} shows that while the TPFA results are strongly influenced by the mesh size, similarly to the finite-volume results from \cite{Raduetal2011}, the unbiased GRW algorithm is practically unconditionally-free of numerical diffusion. The BGRW algorithm is also free of numerical diffusion provided that $\text{P\'{e}}\le 2$ (see also Remark~\ref{remPe}). We also note that $\Delta x=0.05$ defines the coarsest grid acceptable for solving the benchmark problem for coupled flow and transport with BGRW and TPFA codes.

\begin{table}[!ht]
\centering
\caption{Estimation of numerical diffusion for\\ BGRW, GRW and TPFA codes.}
\label{table_relative_errorsNumDiff}
\begin{tabular}{ c c c c c c c  c c c}
\hline
   & $\Delta x$ & $T/\Delta t$ & P\'{e}  & $\varepsilon_{D_{x}}$ & $\varepsilon_{D_{z}}$ \\
  \hline
  \multirow{3}{2.5em}{BGRW}
   &  0.1   & 2   & 3.31 & 7.55e-02 & 2.60e-01 &\\
   &  0.05  & 9   & 1.65 & 1.90e-16 & 1.48e-15 &\\
   &  0.01  & 239 & 0.33 & 4.16e-16 & 1.02e-15 &\\
   &  0.005 & 960 & 0.17 & 2.93e-15 & 3.63e-15 &\\
  \hline
  \multirow{3}{2.5em}{GRW}
   & 0.1    & 4   & 3.31 & 1.94e-16 & 6.14e-16 &\\
   & 0.05   & 4   & 1.65 & 6.60e-17 & 8.05e-16 &\\
   & 0.01   & 19  & 0.33 & 1.94e-16 & 4.79e-16 &\\
   & 0.005  & 39  & 0.17 & 2.10e-15 & 8.92e-16 &\\
  \hline
  \multirow{3}{2.5em}{TPFA}
   &  0.1   & 5   & 3.31 & 9.16e-03 & 1.99e-01 &\\
   &  0.05  & 10  & 1.65 & 4.69e-03 & 9.94e-02 &\\
   &  0.01  & 50  & 0.33 & 9.58e-04 & 1.99e-02 &\\
   &  0.005 & 100 & 0.17 & 5.38e-04 & 9.89e-03 &\\
\hline
\end{tabular}
\end{table}

\subsubsection{Fully coupled water flow and surfactant transport}
\label{coupledFT}

In the following we solve the coupled flow and transport problem (\ref{eqFT01}-\ref{eqFT02}) by using the setup of the benchmark flow problem problem from Section~\ref{2DvalidF} completed by parameters and initial/boundary conditions modeling a situation of coupled water flow and surfactant transport. The surfactant concentration in the domain $\Omega$ has a stratified distribution described by the plane $c(x,z,0)=z/1.2$. Further, the concentration is set to $c=1$ on the Dirichlet boundary $\Gamma_{D_{1}}$ and to $c=0$ on $\Gamma_{D_{2}}$, and no-flow Neumann conditions are imposed on the vertical boundaries.

The flow and transport are coupled in both directions through the van Genuchten-Mualem parameterization (\ref{theta}-\ref{K}) with $\theta(\psi,c)=\theta(\gamma(c)\psi)$, where $\gamma(c)=1/[1-b\ln(c/a+1)]$ models the concentration-dependent surface tension between water
and air \cite{Knabneretal2003}. The constant parameters of $\gamma(c)$ are set to $a=0.44$ and $b=0.0046$ \cite{Illiano2020}. To describe a more realistic heterogeneous soil, the saturated conductivity $K_{sat}$ is modeled as a log-normal space random function with a small variance $\sigma^2=0.5$ and Gaussian correlation of correlation lengths $\lambda_x=0.1$~m and $\lambda_z=0.01$~m in horizontal and vertical directions, respectively. The $\ln K$ field is generated by summing up 100 random periodic modes with the Kraichnan algorithm presented in \cite[Appendix C.3.1.2]{Suciu2019}. The diffusion coefficient is set to a constant value, $D=10^{-3}$~m/day, which is representative for soils and aquifers \cite{Raduetal2011,Schneid2000,Suciu2019}. Following \cite{Illiano2020}, the nonlinear reaction term is specified as $R(c)=10^{-3}c/(1+c)$. Instead of using a fixed number of time steps, as in the flow benchmark presented in Section~\ref{2DvalidF}, now we fix the total time to $T=3$ days, set the intermediate time controlling the drainage process to $\Delta t_D=T/3$, and keep the original time steps $\Delta t$ which ensure the appropriate resolution for contrasting fast and slow processes in loam and clay soils, respectively.

Preliminary tests showed that, in order to obtain an acceptable resolution of the velocity components in the benchmark setup, the unbiased GRW requires extremely fine discretizations with $\Delta x=\mathcal{O}(10^{-5})$. Therefore the transport step is solved with the BGRW algorithm for the mesh size $\Delta x=0.05$ suggested by the above investigations on numerical diffusion. The velocity $\mathbf{q}|_{\partial\Omega}$ on boundaries is approximated by forward finite differences.

\begin{figure}
\begin{minipage}[t]{0.45\linewidth}\centering
\includegraphics[width=\linewidth]{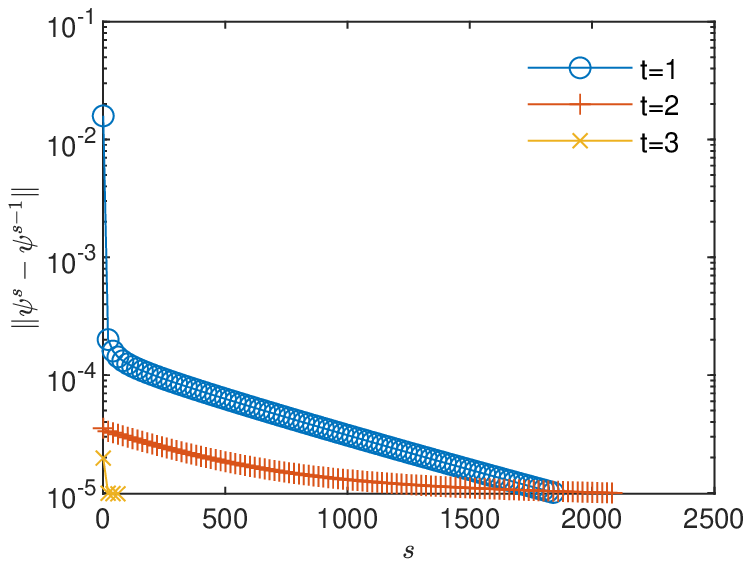}
\caption{\label{figFT11}Convergence of the $L$-scheme implementation of the GRW flow solver for the loam soil problem at three time levels (in days).}
\end{minipage}
\hspace*{0.1in}
\begin{minipage}[t]{0.45\linewidth}\centering
\includegraphics[width=\linewidth]{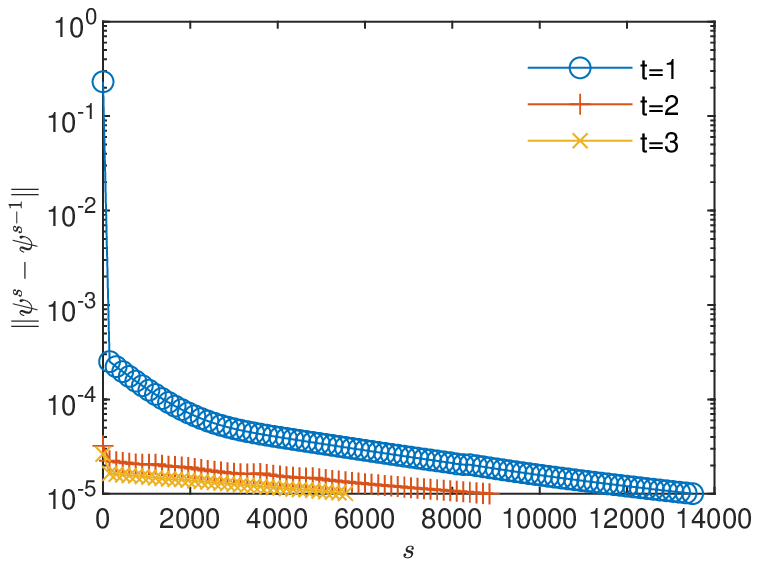}
\caption{\label{figFT12}The same as in Fig.~\ref{figFT11} for the clay soil problem.}
\end{minipage}
\end{figure}

\begin{figure}
\begin{minipage}[t]{0.45\linewidth}\centering
\includegraphics[width=\linewidth]{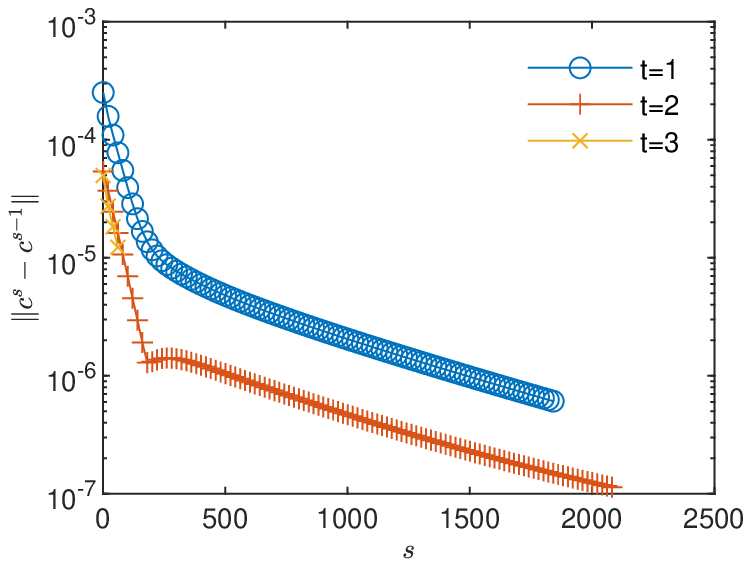}
\caption{\label{figFT13}Convergence of the $L$-scheme implementation of the GRW transport solver for the loam soil problem  at three time levels (in days).}
\end{minipage}
\hspace*{0.1in}
\begin{minipage}[t]{0.45\linewidth}\centering
\includegraphics[width=\linewidth]{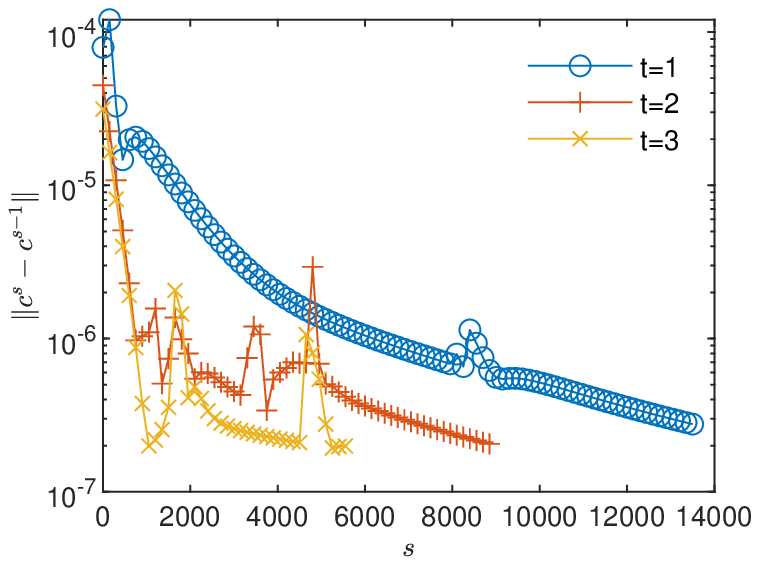}
\caption{\label{figFT14}The same as in Fig.~\ref{figFT13} for the clay soil problem.}
\end{minipage}
\end{figure}

The convergence of the flow and transport $L$-schemes using GRW algorithms requires relatively large linearization parameters, $L_p=L_c=20$, for loam soil, and $L_p=L_c=100$ for clay soil models. These are two order of magnitude larger than for the decoupled-flow benchmark presented in Section~\ref{2DvalidF}, probably due to the increased complexity of the coupled problem. By setting the tolerances of the convergence criterion (\ref{eq5}) to $\varepsilon_a=\varepsilon_r=5\cdot10^{-6}$ the convergence is achieved after about 2000 iterations for the loam soil and about 14000 iterations for the clay soil (see Figs.~\ref{figFT11}~-~\ref{figFT14}).
Estimations of computational orders of convergence \cite{Catinas2019,Catinas2020} indicate slow, power law convergence for both pressure and transport solvers and for both soil models (see Appendix~\ref{appA}).

\begin{figure}
\begin{minipage}[t]{0.45\linewidth}\centering
\includegraphics[width=\linewidth]{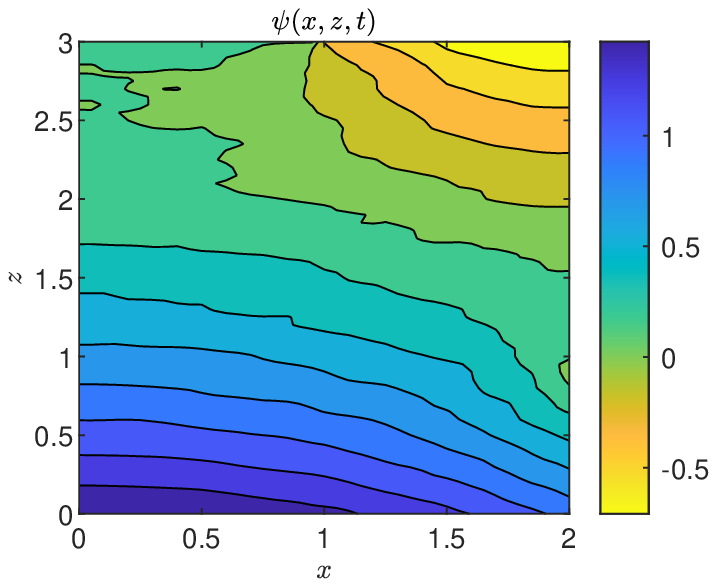}
\caption{\label{figFT1}Pressure head solution $\psi(x,z)$ at $t=T$ for the benchmark problem of recharge from a drainage trench through a silt loam soil coupled with reactive transport.}
\end{minipage}
\hspace*{0.1in}
\begin{minipage}[t]{0.45\linewidth}\centering
\includegraphics[width=\linewidth]{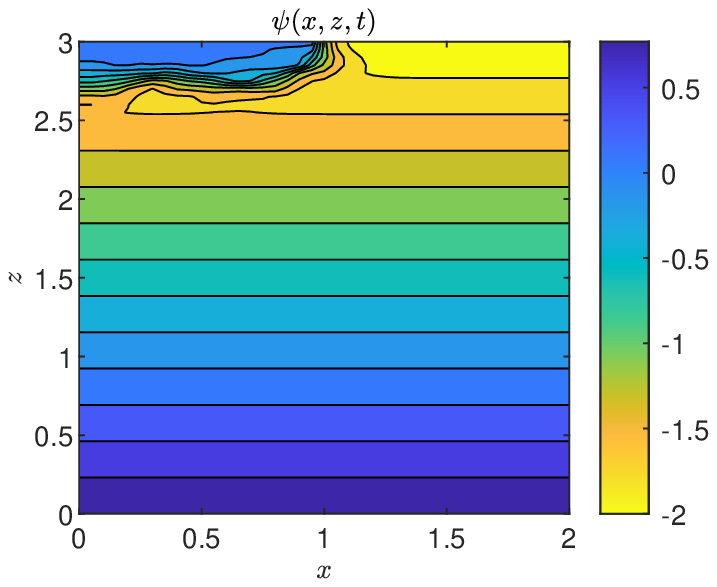}
\caption{\label{figFT2}The same as in Fig.~\ref{figFT1} for a Beit Netofa clay soil.}
\end{minipage}
\end{figure}

\begin{figure}
\begin{minipage}[t]{0.45\linewidth}\centering
\includegraphics[width=\linewidth]{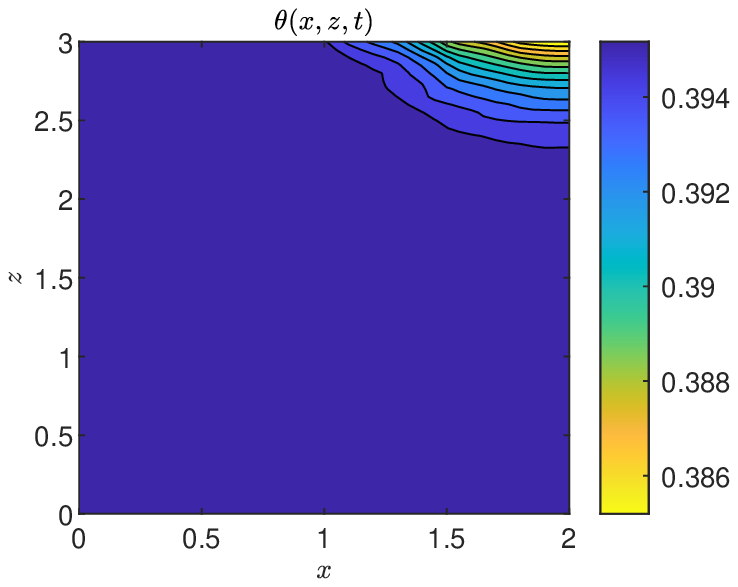}
\caption{\label{figFT3}Water content solution $\theta(x,z)$ at $t=T$ for the benchmark problem of recharge from a drainage trench through a silt loam soil coupled with reactive transport.}
\end{minipage}
\hspace*{0.1in}
\begin{minipage}[t]{0.45\linewidth}\centering
\includegraphics[width=\linewidth]{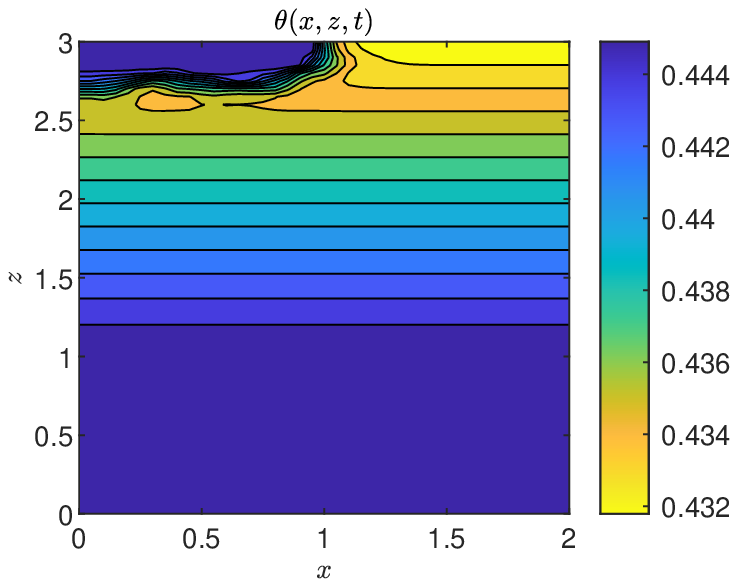}
\caption{\label{figFT4}The same as in Fig.~\ref{figFT3} for a Beit Netofa clay soil.}
\end{minipage}
\end{figure}

\begin{figure}
\begin{minipage}[t]{0.45\linewidth}\centering
\includegraphics[width=\linewidth]{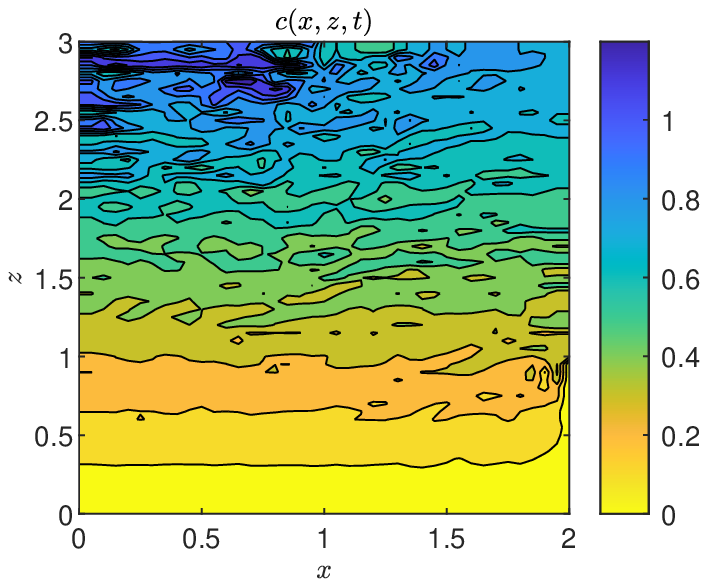}
\caption{\label{figFT5}Concentration solution $c(x,z)$ at $t=T$ for the benchmark problem of recharge from a drainage trench through a silt loam soil coupled with reactive transport.}
\end{minipage}
\hspace*{0.1in}
\begin{minipage}[t]{0.45\linewidth}\centering
\includegraphics[width=\linewidth]{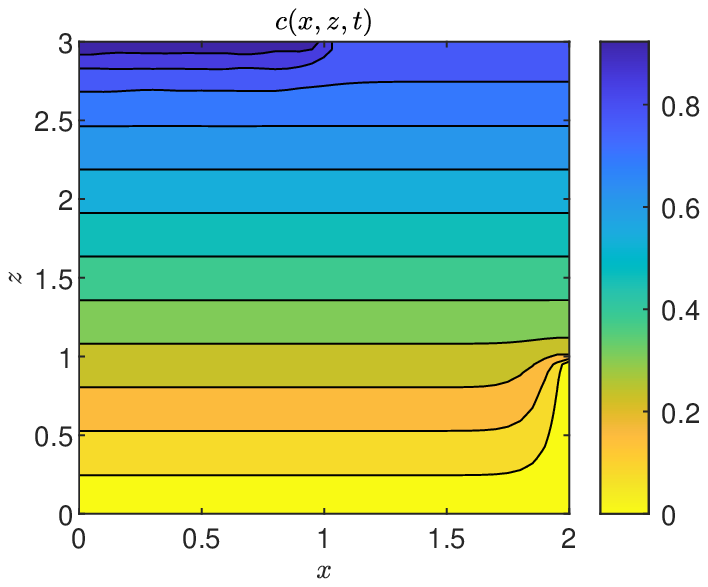}
\caption{\label{figFT6}The same as in Fig.~\ref{figFT5} for a Beit Netofa clay soil.}
\end{minipage}
\end{figure}

\begin{figure}
\begin{minipage}[t]{0.45\linewidth}\centering
\includegraphics[width=\linewidth]{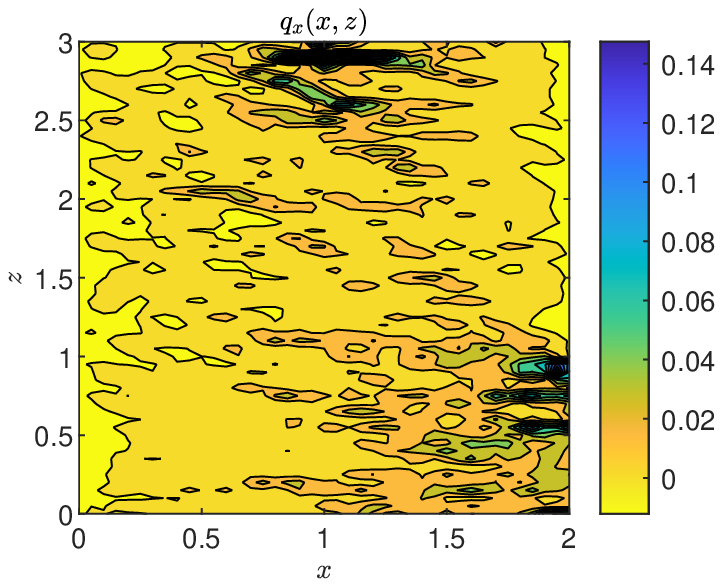}
\caption{\label{figFT7}Horizontal water flux $q_{x}(x,z)$ at $t=T$ for the benchmark problem of recharge from a drainage trench through a silt loam soil coupled with reactive transport.}
\end{minipage}
\hspace*{0.1in}
\begin{minipage}[t]{0.45\linewidth}\centering
\includegraphics[width=\linewidth]{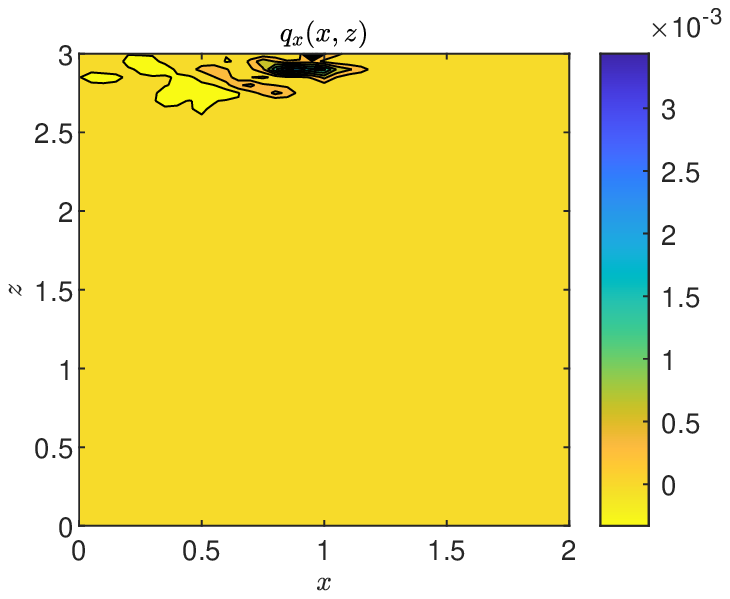}
\caption{\label{figFT8}The same as in Fig.~\ref{figFT7} for a Beit Netofa clay soil.}
\end{minipage}
\end{figure}

\begin{figure}
\begin{minipage}[t]{0.45\linewidth}\centering
\includegraphics[width=\linewidth]{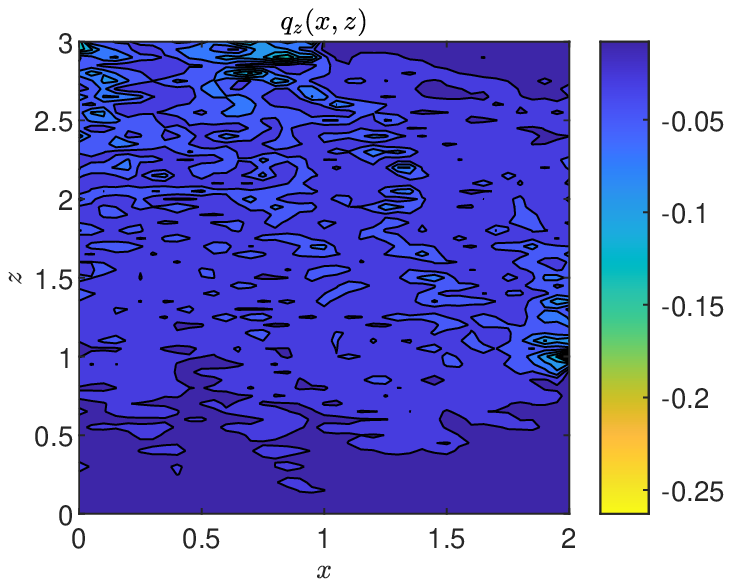}
\caption{\label{figFT9}Vertical water flux $q_{x}(x,z)$ at $t=T$ for the benchmark problem of recharge from a drainage trench through a silt loam soil coupled with reactive transport.}
\end{minipage}
\hspace*{0.1in}
\begin{minipage}[t]{0.45\linewidth}\centering
\includegraphics[width=\linewidth]{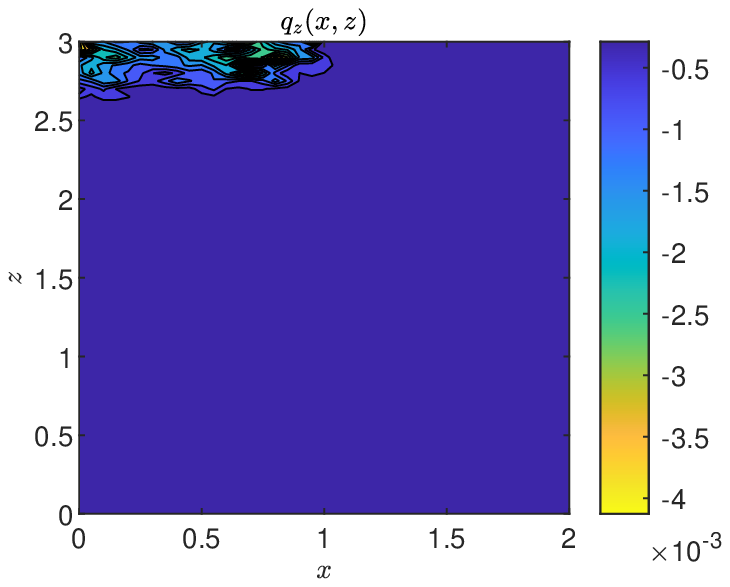}
\caption{\label{figFT10}The same as in Fig.~\ref{figFT9} for a Beit Netofa clay soil.}
\end{minipage}
\end{figure}

The results obtained by coupling the GRW-flow and BGRW-transport solvers are presented in Figs.~\ref{figFT1}-\ref{figFT10}. The randomness of $K_{sat}$ is especially felt by the pressure distribution in the more permeable loam soil (Fig.~\ref{figFT1}), while in the clay soil the pressure remains almost stratified (Fig.~\ref{figFT2}). The same contrast is shown by the water content, with almost saturated loam soil (Fig.~\ref{figFT3}) and partially stratified saturation in the clay soil (Fig.~\ref{figFT4}). Since the Darcy velocity is proportional to the gradient of the random pressure, the heterogeneity of the advective component of the transport process is mainly manifest in the final distribution of the concentration in the loam and clay soils (compare Figs.~\ref{figFT5} and Fig.~\ref{figFT6}). Significant differences between the loan and clay soils are also illustrated by the spatial distribution of the velocity components (Figs.~\ref{figFT7}~-~\ref{figFT10}).

The results obtained with the GRW/BGRW flow and transport solvers are compared with those provided by a TPFA code using $L_p=L_c=1$, for both soils, and $L_c=2L_p$. The convergence is achieved in reasonable computing times of 263 seconds (loam) and 177 seconds (clay) only when using the Anderson acceleration procedure \cite{Anderson1965,WalkerandNi2011,Both2019}. Note that the GRW times on the same computer are of the same order of magnitude (526 and 178, respectively), without appealing to the acceleration procedure.

The errors for pressure, water content and velocity components shown in Table~\ref{table_relative_errorsFT} are more or less similar to those for the flow benchmark problem given in Table~\ref{table_relative_errors}. The difference of one order of magnitude between the $\varepsilon_c$ values for the two soils can be traced back to the amount of numerical diffusion of the TPFA transport solver (see Table~\ref{table_relative_errorsNumDiff}). The estimated mean P\'{e}clet number for the loam soil, $\text{P\'{e}}\approx 1.3$, is much larger than the value $\text{P\'{e}}\approx 4\cdot 10^{-3}$ estimated for the clay soil and can partially explain the larger $\varepsilon_c$ value in the first case. Since the pressure equation is essentially an advection-diffusion equation with velocity given by the derivatives of the coefficient $K$ (see e.g., \cite{Gotovacetal2009,Suciu2020}), the errors $\varepsilon_{q_x}$ and $e_{q_z}$, of order $~10^{-1}$ also could be produced by numerical diffusion, in the flow solver. In the setup of the benchmark problems, for both coupled flow and transport and decoupled flow, we estimate a mean P\'{e}clet number $\text{P\'{e}}\approx 0.9$ for both loam and clay soil models (for comparison, in the one dimensional case with smaller $\varepsilon_q$,  $\text{P\'{e}}$ was about 0.03 in Scenario 1 and 0.3 in Scenario 2). Since the flow and transport solvers implemented in MRST basically use the same TPFA finite volume method, we may expect that the flow solver produces a numerical diffusion comparable to that of the transport solver shown in Table~\ref{table_relative_errorsNumDiff}.

\begin{table}[!ht]
\centering
\caption{Comparison of GRW and TPFA solutions\\ of the coupled
flow-transport benchmark problem.}
\label{table_relative_errorsFT}
\begin{tabular}{ c c c c c c  c c c}
\hline
& $\varepsilon_{\psi}$  & $\varepsilon_c$ & $\varepsilon_\theta$  & $\varepsilon_{q_x}$ & $\varepsilon_{q_z}$\\ \hline
loam   & 2.89e-02 & 4.79e-01  & 7.25e-05 & 3.15e-01  & 2.18e-01\\
clay   & 5.95e-02 & 3.77e-02  & 7.61e-04 & 3.66e-01  & 5.36e-01\\ \hline
\end{tabular}
\end{table}

\subsubsection{Comparison with one-dimensional coupled flow and transport solutions}
\label{1Dvalid_flow_transp}

Solutions to one-dimensional problems are readily obtained by particularizing the algorithms and the methodology used in the two-dimensional case. For code verification purposes, we consider the manufactured solutions
\begin{equation}\label{eqFTsol_1D}
\psi_m(z,t) = - t\ z\ (z-1) -\ 1, \quad
c_m(z,t) = t\ z\ (z-1) +\ 1,
\end{equation}
and the parameter functions
\begin{equation}\label{eqFTparam_1D}
\theta(\psi,c) = \frac{1}{1-\psi - c/10}\;, \quad K(\theta(\psi)) = \psi^{2}\;.
\end{equation}

The coupled system of equations (\ref{eqFT01}-\ref{eqFT02}) is solved in the unit interval $[0,1]$ for a total time $T=1$. The initial and boundary conditions, as well as the source term, are specified by the exact solutions (\ref{eqFTsol_1D}). With tolerances set to $\varepsilon_a=10^{-6}$ and $\varepsilon_r=0$, the convergence criterion (\ref{eq5}) is satisfied for $L_p=L_c=50$ after about 700 iterations.

\begin{table}[h!]
\begin{center}
\caption{Estimated order of convergence of the one-dimensional GRW solver: pressure solutions.}
\label{tab:mrst-grwFTp1D}
\begin{tabular}{c c c c c c c c}
\hline
&$\varepsilon_1$ & EOC & $\varepsilon_2$ & EOC & $\varepsilon_3$ & EOC & $\varepsilon_4$\\
\hline
GRW (analytical $\mathbf{q}|_{\partial\Omega}$) & 3.41e-02 & 2.21 & 7.35e-03 & 2.00 & 1.83e-03 & 1.64 & 5.90e-04 \\
\hline
GRW (approximate $\mathbf{q}|_{\partial\Omega}$) & 3.41e-02 & 2.21 & 7.34e-03 & 2.01 & 1.82e-03 & 1.66 & 5.77e-04 \\
\hline
GRW ($\mathbf{q}|_{\partial\Omega}$ from $\mbox{int}(\Omega)$) & 3.41e-02 & 2.22 & 7.34e-03 & 2.01 & 1.82e-03 & 1.66 & 5.76e-04 \\
\hline 
\end{tabular}
\end{center}
\end{table}

\begin{table}[h!]
\begin{center}
\caption{Estimated order of convergence of the one-dimensional GRW solver: concentration solutions.}
\label{tab:mrst-grwFTc1D}
\begin{tabular}{c c c c c c c c}
\hline
&$\varepsilon_1$ & EOC & $\varepsilon_2$ & EOC & $\varepsilon_3$ & EOC & $\varepsilon_4$\\
\hline
GRW (analytical $\mathbf{q}|_{\partial\Omega}$) & 2.10e-02 & 2.12 & 4.85e-03 & 2.12 & 1.12e-03 & 2.83 & 1.56e-04 \\
\hline
GRW (approximate $\mathbf{q}|_{\partial\Omega}$) & 4.24e-02 & 2.08 & 1.00e-02 & 1.94 & 2.61e-03 & 1.25 & 1.10e-03 \\ \hline
GRW ($\mathbf{q}|_{\partial\Omega}$ from $\mbox{int}(\Omega)$) & 4.06e-02 & 2.20 & 8.85e-03 & 1.98 & 2.25e-03 & 1.17 & 1.00e-03 \\
\hline
\end{tabular}
\end{center}
\end{table}

Comparing with the results of code verification tests in the two-dimensional case presented in Tables~\ref{tab:mrst-grwFTp}~and~\ref{tab:mrst-grwFTc}, one remarks the same tendency to convergence of order 2 of the EOC values in one-dimensional case (Tables~\ref{tab:mrst-grwFTp1D}~and~\ref{tab:mrst-grwFTc1D}). Instead, the successive absolute errors $\varepsilon^{(k)}$ with respect to the manufactured solutions $\psi_m$ and $c_m$ are systematically larger by one order of magnitude than in the two-dimensional case.

Additionally, we test the ability of the one-dimensional GRW solvers to describe the transition from unsaturated to saturated regime. Following \cite{Raduetal2009}, we set $K=1$ and solve the degenerate Richards equation coupled with the transport equation from the previous example. We maintain the bidirectional coupling by choosing the water content as a function of both $\psi$ and $c$, $\theta=1/(3.4333-p-c/10)$ for $\psi<0$ and $\theta=0.3$ for $\psi\ge 0$. We keep the analytical solution $c_m$ for concentration unchanged and chose the pressure solution as $\psi_m=-tx(1-x)+x/4$, such that $\psi<0$ for $x\in(0,1-1/(4t))$ and $\psi\ge 0$ otherwise. Since the fully coupled flow and transport problem is now degenerate, the source functions obtained after inserting $\psi_m$ and $c_m$ into the equations (\ref{eqFT01}-\ref{eqFT02}) also possess two branches, corresponding to $\psi<0$ and $\psi\ge 0$. For the same tolerances as in the previous example, $\varepsilon_a=10^{-6}$, $\varepsilon_r=0$, and $L_p=L_c=100$ the convergence criterion (\ref{eq5}) is satisfied after a number of iterations ranging between 400 and 2780. The results presented in Tables~\ref{tab:mrst-grwFTp1D_deg}~and~\ref{tab:mrst-grwFTc1D_deg} demonstrate the convergence of the GRW solvers in degenerate case and indicate the tendency to the 2-nd order of convergence.

\begin{table}[h!]
\begin{center}
\caption{Estimated order of convergence of the one-dimensional GRW solver for a degenerate problem: pressure solutions.}
\label{tab:mrst-grwFTp1D_deg}
\begin{tabular}{c c c c c c c c}
\hline
&$\varepsilon_1$ & EOC & $\varepsilon_2$ & EOC & $\varepsilon_3$ & EOC & $\varepsilon_4$\\
\hline
GRW (analytical $\mathbf{q}|_{\partial\Omega}$) & 4.59e-02 & 2.00 & 1.14e-02 & 1.95 & 2.95e-03 & 1.57 & 9.97e-04 \\
\hline
GRW (approximate $\mathbf{q}|_{\partial\Omega}$) & 4.59e-02 & 2.01 & 1.14e-02 & 1.96 & 2.94e-03 & 1.56 & 9.99e-04 \\
\hline
GRW ($\mathbf{q}|_{\partial\Omega}$ from $\mbox{int}(\Omega)$) & 4.59e-02 & 2.01 & 1.14e-02 & 1.96 & 2.94e-03 & 1.56 & 9.99e-04 \\
\hline 
\end{tabular}
\end{center}
\end{table}

\begin{table}[h!]
\begin{center}
\caption{Estimated order of convergence of the one-dimensional GRW solver for a degenerate problem: concentration solutions.}
\label{tab:mrst-grwFTc1D_deg}
\begin{tabular}{c c c c c c c c}
\hline
&$\varepsilon_1$ & EOC & $\varepsilon_2$ & EOC & $\varepsilon_3$ & EOC & $\varepsilon_4$\\
\hline
GRW (analytical $\mathbf{q}|_{\partial\Omega}$) & 7.64e-03 & 1.97 & 1.95e-03 & 2.16 & 4.38e-04 & 2.28 & 8.99e-05 \\
\hline
GRW (approximate $\mathbf{q}|_{\partial\Omega}$) & 5.02e-02 & 1.94 & 1.31e-02 & 1.90 & 3.51e-03 & 1.37 & 1.36e-03 \\ \hline
GRW ($\mathbf{q}|_{\partial\Omega}$ from $\mbox{int}(\Omega)$) & 5.38e-02 & 1.91 & 1.43e-02 & 1.90 & 3.84e-03 & 1.42 & 1.44e-03 \\
\hline
\end{tabular}
\end{center}
\end{table}

Similarly to the benchmark problem from Section~\ref{coupledFT} and using the same physical parameters, a one one-dimensional problem of coupled flow and transport is solved for $z\in(0,3)$, with initial conditions $p(z,0)=1-z$, $c(z,0)=z/1.2$, boundary condition $p(3,t)$ given by that imposed on $\Gamma_{D_{1}}$ in two dimensional case, and $p(0,t)=p(0,0)=1$ (water table fixed at $z=1$, as in the two-dimensional case).

\begin{figure}[h]
\begin{minipage}[t]{0.45\linewidth}\centering
\includegraphics[width=\linewidth]{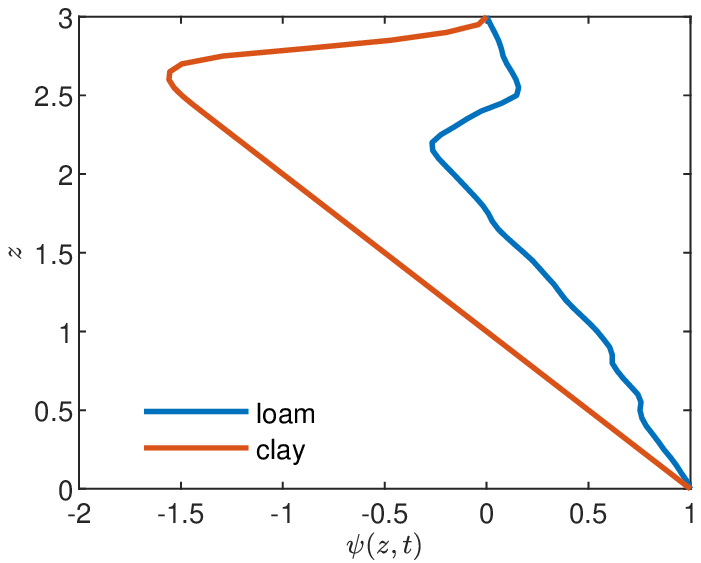}
\caption{\label{fig_p1D}Pressure head at $t=T$ in loam and clay soils.}
\end{minipage}
\hspace*{0.1in}
\begin{minipage}[t]{0.45\linewidth}\centering
\includegraphics[width=\linewidth]{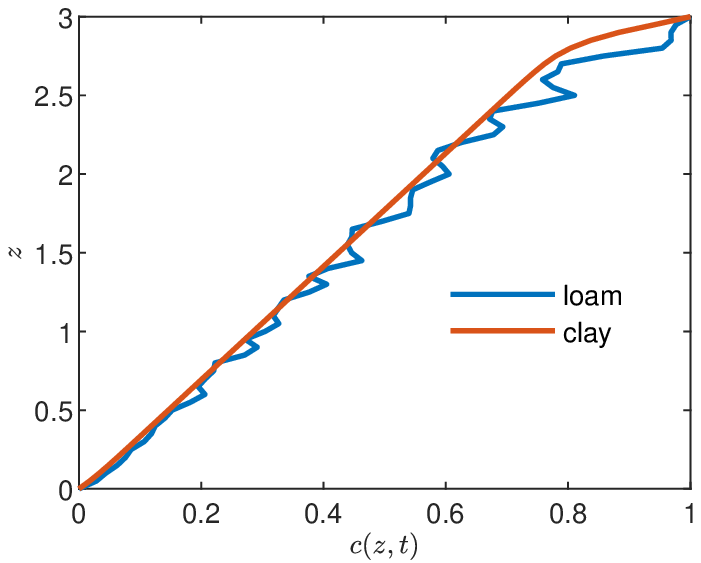}
\caption{\label{fig_c1D}Concentration at $t=T$ in loam and clay soils.}
\end{minipage}
\end{figure}

\begin{figure}[h]
\begin{minipage}[t]{0.45\linewidth}\centering
\includegraphics[width=\linewidth]{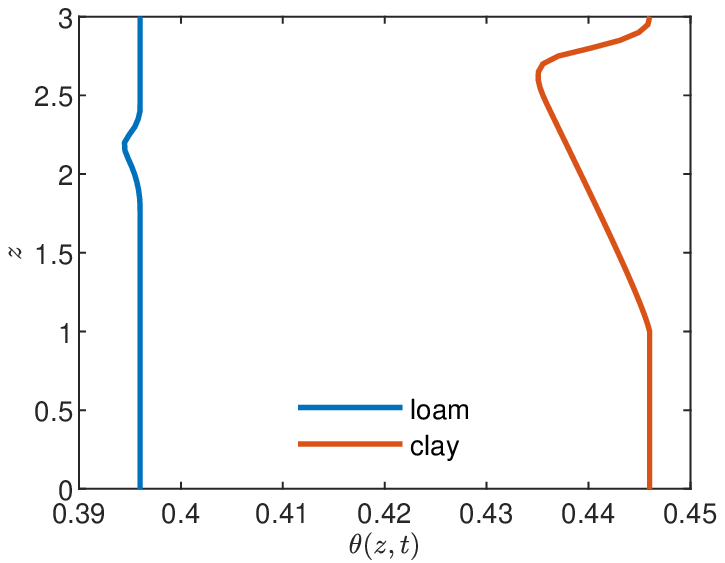}
\caption{\label{fig_t1D}Water content at $t=T$ in loam and clay soils.}
\end{minipage}
\hspace*{0.1in}
\begin{minipage}[t]{0.45\linewidth}\centering
\includegraphics[width=\linewidth]{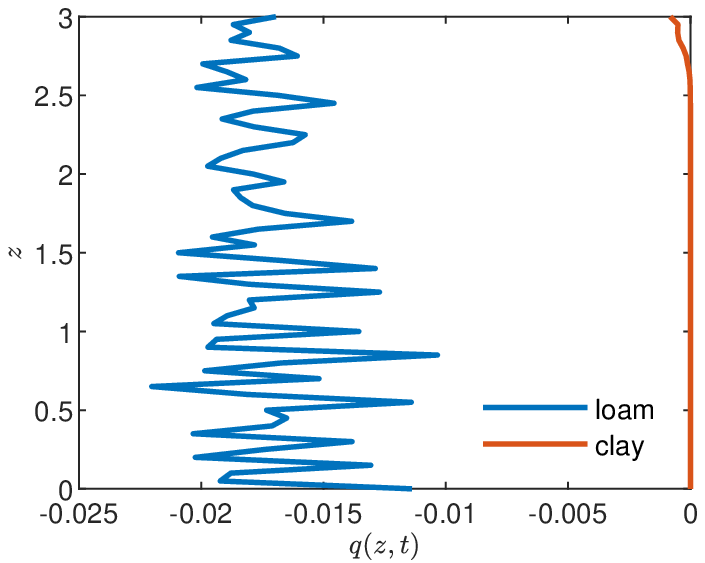}
\caption{\label{fig_q1D}Water flux at $t=T$ in loam and clay soils.}
\end{minipage}
\end{figure}

The solutions for pressure head, concentration, water content and water flux for loam and clay soils are compared in Figs.~\ref{fig_p1D}~-~\ref{fig_q1D}. One remarks that, even though the lateral heterogeneity of the two-dimensional benchmark is ignored, the main features are also revealed  by the one-dimensional drainage model: the discrepancy between fast-loam and slow-clay flow and transport processes, the same intervals of variation of the solutions, and similar behavior on the vertical direction.

\section{Two-dimensional GRW solutions for groundwater flow and transport at regional and field scales}
\label{2Dreg_scale}

For saturated aquifers ($\theta=const$) Eq.~(\ref{eq11}) reduces to a linear equation solved by the steady state hydraulic head solution in $h(x,y)$, under time independent boundary conditions. As noted in Remark~\ref{remFsat2D}, the GRW $L$-scheme (\ref{eq12}-\ref{eq14}) becomes, in this case, a transient scheme for the linear flow equation. In the following examples, we consider flow problems formulated in two-dimensional domains, $(x,y)\in[0,L_x]\times[0,L_y]$, with Dirichlet boundary conditions $h(0,y)=H_1$ and $h(L_x,y)=H_2$ and no-flow Neumann conditions on top and bottom boundaries. In the saturated flow regime, the transport Eq.~(\ref{eqFT02}) is also linear and decoupled from the linear flow equation. Decoupled transport problems can be solved by either biased- or unbiased-GRW algorithms (see Remark~\ref{remTdec} and Section~\ref{2DalgFT_GRW}) on the same lattice as that used to compute the flow velocity.

\subsection{Flow in heterogeneous aquifers at regional scale}
\label{2Dreg_scale_heterogeneous}

For the beginning, we follow the setup for regional scale used in \cite{Franssenetal2009} to compare approaches for inverse modeling of groundwater flow. The domain and the boundary conditions are specified by $L_x=4900$~m, $L_y=5000$~m, $H_1=0$~m, $H_2=5$~m. The hydraulic conductivity $K$ is a log-normally distributed random field defined by the mean $\langle K\rangle=12\cdot 10^{-4}$~m/s, the correlation length $\lambda=500$ m, and the variance $\sigma^2=1$ of the $\ln K$-field. The $K$-field is generated, as in Section~\ref{coupledFT} above, by summing 100 random periodic modes with the Kraichnan algorithm. Besides the exponential correlation considered in \cite{Franssenetal2009}, we also investigate the behavior of the flow solution for Gaussian correlation of the $\ln K$ field with the same correlation length, as well as in case of the smaller variance $\sigma^2=0.1$, for both correlation models.

The two correlation models of the $\ln K$-field are of the form $C(r)=\sigma^2\exp[-(r/\lambda)^\alpha]$,
where $r=(r_{x}^{2}+r_{y}^{2})^{1/2}$ is the spatial lag, the exponent $\alpha=1$ corresponds to the exponential model, and $\alpha=2$ to the Gaussian one. Since the correlation functions depend on spatial variables through $r/\lambda$, the computation can be done for spatial dimensions scaled by $\lambda$, that is, fields of dimensionless correlation length $\lambda^*=1$ and a domain $[0,L_x/\lambda]\times[0,L_y/\lambda]$. The results on the original grid are finally obtained after the multiplication by $\lambda$ of the solution $h(x,y)$ and of the spatial coordinates.

\begin{figure}[h]
\begin{minipage}[t]{0.45\linewidth}\centering
\includegraphics[width=\linewidth]{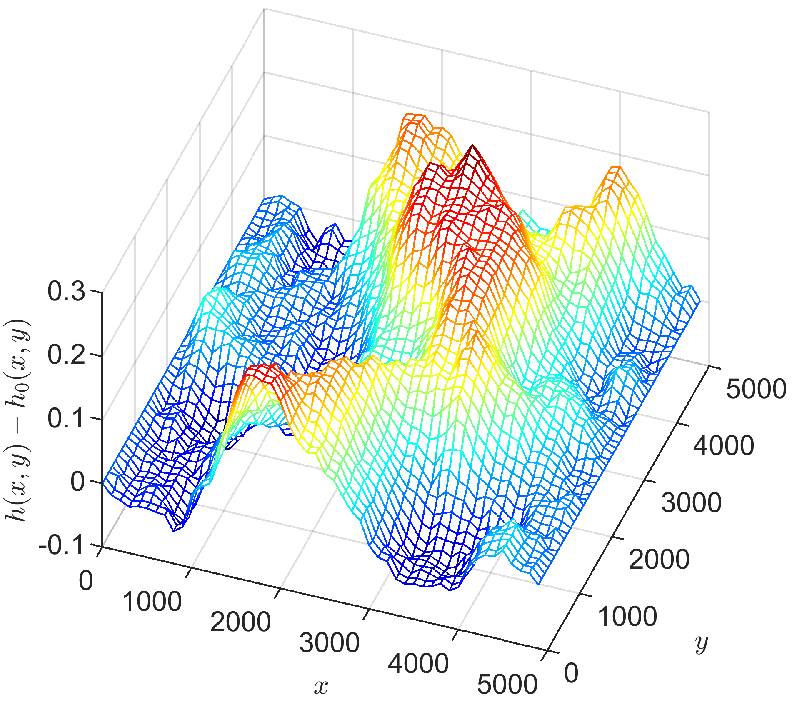}
\caption{\label{fig:h_regional}Fluctuations of the hydraulic head $h(x,y)$ solution obtained with scaled spatial dimensions computed for a fixed realization of the $\ln(K)$ field with exponential correlation.}
\end{minipage}
\hspace*{0.1in}
\begin{minipage}[t]{0.45\linewidth}\centering
\includegraphics[width=\linewidth]{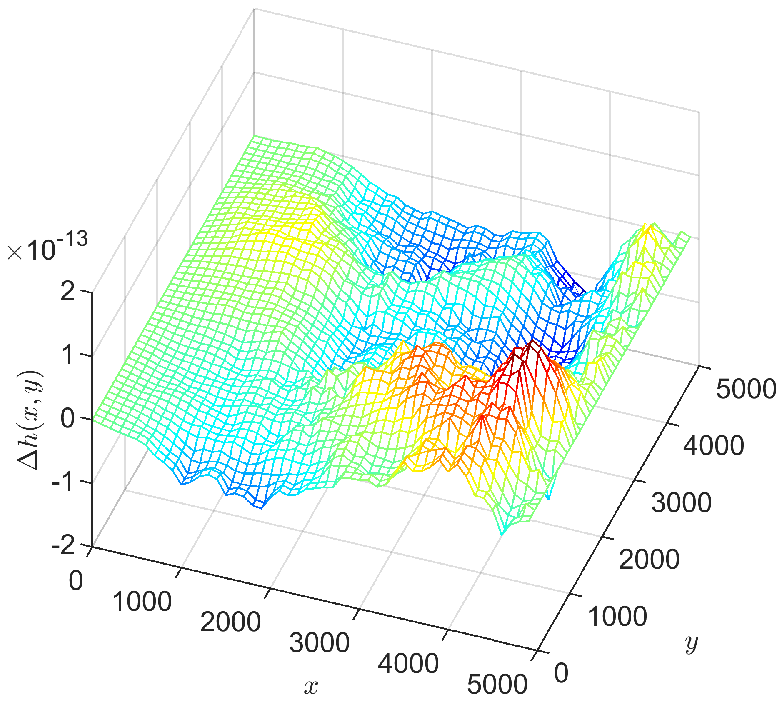}
\caption{\label{fig:dh_regional} Deviation $\Delta h(x,y)=h(x,y)-\hat{h}(x,y)$ of the hydraulic head solution of the scaled problem from the solution $\hat{h}(x,y)$ of the original problem.}
\end{minipage}
\end{figure}

The solutions $h(x,y)$ of the stationary equation~(\ref{eq11}) corresponding to $\theta=const$, for given realizations of the $K$-field with $\sigma^2=0.1$, are obtained under the initial condition $h_0(x,y)$, which is the plane defined by the Dirichlet boundary conditions $h(0,y)=0$ and $h(L_x/\lambda,y)=H_2/\lambda$. To test the setup for the solution obtained with the geometry scaled by $\lambda$ (computed with $\Delta x=\Delta y=0.2$~m), shown in Fig.~\ref{fig:h_regional}, we compare the solution $h$ to the unscaled solution $\hat{h}$ for the same realization of the $K$-field with exponential correlation and $\sigma^2=0.1$ (computed with $\Delta x=\Delta y=100$~m). The steady state is reached after about $4\cdot 10^{5}$ iterations of the GRW solver, for both scaled and unscaled geometry. The difference of the two solutions shown in Fig.~\ref{fig:dh_regional} and the corresponding relative error $\|\Delta h\|/\|\hat{h}\|=1.29\cdot 10^{-14}$ are close to the machine precision.

To estimate the order of convergence of the GRW scheme for this particular flow problem, we use manufactured analytical solutions provided in the Git repository \href{https://github.com/PMFlow/FlowBenchmark}{https://github.com/PMFlow/FlowBenchmark} and, similarly to estimations performed in Section~\ref{2DvalidF}, we compute the EOC according to (\ref{eq16}) by successively halving the space steps from $\Delta x=\Delta y=2\cdot 10^{-1}$ up to $\Delta x=\Delta y=2.5\cdot 10^{-2}$.

\begin{table}[!ht] \centering
 \caption{Computational order of convergence of the GRW scheme estimated according to (\ref{eq16}).}
 \label{EOC}
 \begin{tabular}{ c c c c c c c c c c c c c }
  \hline
  Correlation model & $\sigma^2$ & $\varepsilon_1$ & EOC & $\varepsilon_2$ & EOC &$\varepsilon_3$ & EOC & $\varepsilon_4$\\
  \hline
  \multirow{2}{2em}{Exponential}
  & $0.1$ & 1.35e+01 & 3.67 & 1.06e+00 & 1.86 & 2.92e-01 & 0.66 & 1.85e-01\\
  & $1$   & 1.80e+02 & 3.24 & 1.90e+01 & 2.09 & 4.47e+00 & 1.96 & 1.15e+00\\
  \hline
  \multirow{2}{2em}{Gaussian}
  & $0.1$ & 7.37e-02 & 1.98 & 1.87e-02 & 1.63 & 6.03e-03 & 1.14 & 2.73e-03\\
  & $1$   & 1.31e-01 & 1.59 & 4.35e-02 & 1.51 & 1.53e-02 & 1.47 & 5.51e-03\\
  \hline
 \end{tabular}
\end{table}

\begin{table}[!ht] \centering
 \caption{Computational order of convergence of the TPFA solver estimated according to (\ref{eq16}).}
 \label{EOC_MRST}
 \begin{tabular}{ c c c c c c c c c c c c c }
  \hline
  Correlation model & $\sigma^2$ & $\varepsilon_1$ & EOC & $\varepsilon_2$ & EOC &$\varepsilon_3$ & EOC & $\varepsilon_4$\\
  \hline
  \multirow{2}{2em}{Exponential}
  & $0.1$ & 4.67e+00 & 1.71 & 1.43e+00 & 1.95 &3.70e-01  & 0.48 & 2.65e-01\\
  & $1$   & 1.01e+02 & 2.23 & 2.14e+01  &3.11 &2.48e+00  & 0.41 & 1.86e+00\\
  \hline
  \multirow{2}{2em}{Gaussian}
  & $0.1$ & 9.22e-02 & 2.00 & 2.30e-02 & 2.00 & 5.75e-03 & 2.00 & 1.44e-03\\
  & $1$   & 1.84e-01 & 2.00 & 4.61e-02 & 2.00 & 1.16e-02 & 2.00 & 2.89e-03\\
  \hline
 \end{tabular}
\end{table}

We note that the EOC approach presented here differers somewhat from that used in \cite{Alecsaetal2020,Suciu2020}. The reference solution is now the manufactured solution, instead of the solution on the finest grid, and the error norm is no longer computed after the first iteration but after large numbers of iterations (from $10^5$ to more than $10^7$ ), when the GRW solution approaches the stationarity. Due to the limited number of iterations, the solutions are not yet strictly stationary and the order of convergence may be not accurately estimated in some cases. Therefore we also use a TPFA flow solver to compute EOC values for the same Scenarios.

The results presented in Tables~\ref{EOC}~and~\ref{EOC_MRST} show significant differences between the two correlation models. For Gaussian correlation the errors obtained with the two approaches are relatively small in all cases. Instead, for exponential correlation, despite the strong EOC obtained after the first two refinements, the errors are extremely large for $\sigma^2=1$ and become smaller than one only for $\sigma^2=0.1$, after the second refinement of the grid. These results are consistent with those presented in \cite{Alecsaetal2020}, where similar benchmark problems were solved for a larger range of parameters of the $\ln K$ field.

\subsection{Flow in conditions of random recharge}
\label{2Dreg_scale_recharge}

We consider in the following a flow problem formulated for the same geometry and boundary conditions as in the previous subsection, which has been used in \cite{Pasettoetal2011} to design a new Monte Carlo approach for flow driven by spatially distributed stochastic sources. Now the hydraulic conductivity is constant, $K=12\cdot 10^{-4}$~m/s, and the groundwater recharge is described by a source term $f$ in Eq.~(\ref{eq11}), modeled as a random space function of mean $\langle f\rangle=362.912$~mm/year, log-normally distributed with exponential correlation specified by different correlation lengths and variances of the $\ln f$ field. Among  different scenarios presented in\cite{Pasettoetal2011}, we consider for comparison with the present computations only the case $\lambda=500$~m and the variance $\sigma^2=1$.

To test the scaling procedure, we also consider a sink term $Q$ located at $(25\lambda,25\lambda)$.
\begin{figure}[h]
\begin{minipage}[t]{0.45\linewidth}\centering
\includegraphics[width=\linewidth]{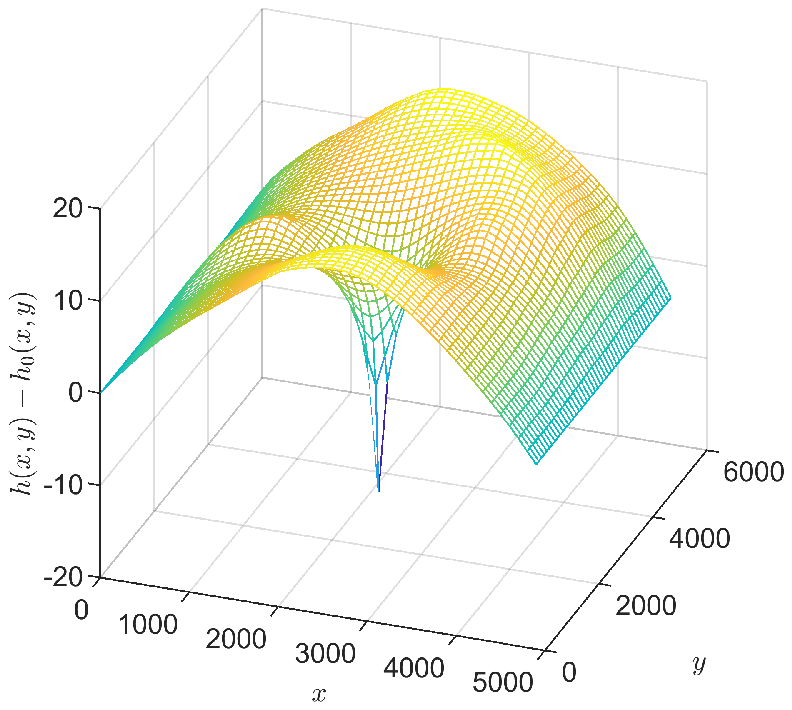}
\caption{\label{fig:h_recharge}Hydraulic head $h(x,y)$ solution of the random recharge problem with scaled spatial dimensions computed with a fixed realization of the $\ln(K)$ field with exponential correlation.}
\end{minipage}
\hspace*{0.1in}
\begin{minipage}[t]{0.45\linewidth}\centering
\includegraphics[width=\linewidth]{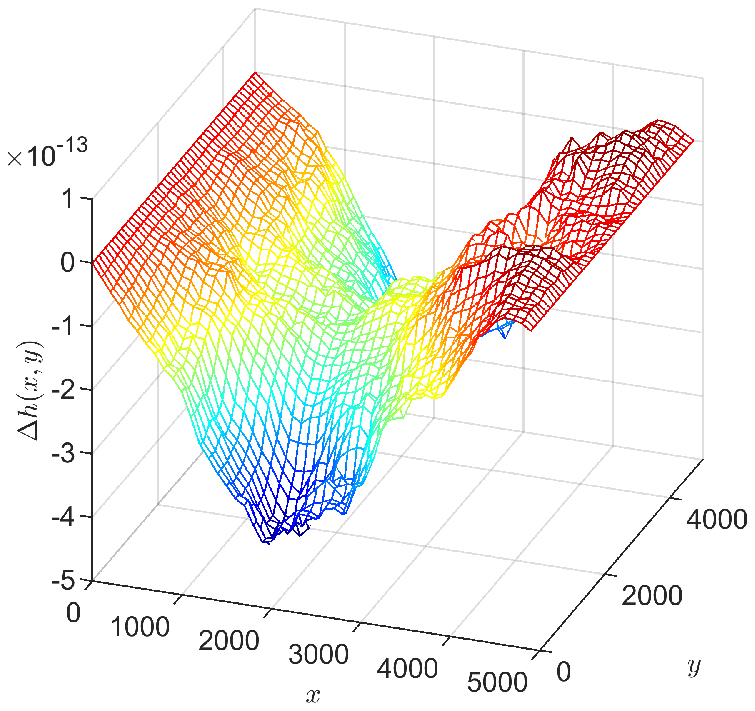}
\caption{\label{fig:dh_recharge} Deviation $\Delta h(x,y)=h(x,y)-\hat{h}(x,y)$ of the hydraulic head solution of the scaled problem from the solution $\hat{h}(x,y)$ of the original problem.}
\end{minipage}
\end{figure}
The setup for the solution obtained with the geometry scaled by $\lambda$, shown in Fig.~\ref{fig:h_recharge}, is evaluated in the same way as in the previous subsection. The difference of the two solutions shown in Fig.~\ref{fig:dh_recharge} and the corresponding relative error $\|\Delta h\|/\|\hat{h}\|=7.14\cdot 10^{-15}$ are again close to the machine precision.

In a first validation test, we compare the GRW and TPFA solutions of the random recharge problem on the computational domain scaled by $\lambda=500$~m, for single-realizations of the random recharge with both exponential and Gaussian correlation of the $\ln f$ field and two variances, $\sigma^2=0.1$ and $\sigma^2=1$. The absolute and relative differences, $\varepsilon_{a}=\| h^{GRW} - h^{TPFA}\|$ and $\varepsilon_{r}=\|h^{GRW} - h^{TPFA}\|/\|h^{TPFA} \|$,  presented in Table~\ref{tab:MRST_GRW} indicate a good agreement between the two approaches.
\begin{table}[h!]
\begin{center}
\caption{Comparison of GRW and TPFA solutions of the\\ random recharge problem.}
\label{tab:MRST_GRW}
\begin{tabular}{c c c c}
\hline
 Correlation model & $\sigma^2$ & $\varepsilon_{a}$ & $\varepsilon_{r}$ \\
\hline
 \multirow{2}{2em}{Exponential}
& 0.1 & 63.44  & 5.97e-2  \\
& 1   & 101.71 & 9.82e-2  \\
\hline
  \multirow{2}{2em}{Gaussian}
& 0.1 &  84.12 & 8.72e-2 \\
& 1   & 137.09 & 1.62e-2 \\
\hline
\end{tabular}
\end{center}
\end{table}

Further, we perform statistical inferences of the mean and variance obtained from an ensemble of 100 Monte Carlo simulations within the setup of \cite{Pasettoetal2011} for random recharge term with exponential correlation and variance $\sigma^2=1$. The mean and the variance of the hydraulic head $h$ are computed as averages over realizations of the $\ln f$ field followed by spatial averages, with standard deviation estimated by spatial averaging. The results presented in Table~\ref{tableMCh} show, again, that the GRW and TPFA results are in good statistical agreement.
\begin{table}[!ht]
\centering
\caption{Statistical moments of the hydraulic head \\
(Monte Carlo and spatial averages).}
\label{tableMCh}
\begin{tabular}{ c c c c c }
\hline
& mean  & variance &\\ \hline
GRW  & 21.51$\pm$9.17 & 41.11$\pm$27.82 &\\ \hline
TPFA & 19.74$\pm$7.84 & 32.09$\pm$21.33 &\\ \hline
\end{tabular}
\end{table}

Finally, we compare the mean and the variance estimated at the center of the computational domain by GRW and TPFA simulations with the results presented in \cite{Pasettoetal2011}. As seen in Table~\ref{table_compPasseto}, the mean values compare quite well but both the GRW and TPFA approaches overestimate the variance computed for the same parameters in \cite[Fig. 6]{Pasettoetal2011}. This discrepancy can be attributed either to the large errors expected for exponential correlation model (see Tables~\ref{EOC} and \ref{EOC_MRST}) or to the statistical inhomogeneity of the Monte Carlo ensemble of 100 realizations indicated by the large standard deviations shown in Table~\ref{tableMCh}.
\begin{table}[!ht]
\centering
\caption{Statistical moments of the hydraulic head \\
(MC averages at the center of the domain).}
\label{table_compPasseto}
\begin{tabular}{ c c c c c }
\hline
& mean  & variance &\\ \hline
GRW                    & 31.67 & 65.14 &\\ \hline
TPFA                   & 28.31 & 53.39 &\\ \hline
(Passeto et al., 2011) & 31.05 & 40.08 &\\ \hline
\end{tabular}
\end{table}

\subsection{Flow and advection-dispersion transport in aquifers}
\label{2D_FlowTransp_aquifers}

In the following we consider an incompressible flow in the domain $[0,20]\times[0,10]$, driven by Dirichlet boundary conditions $h(0,y)=1$ and $h(20,y)=0$ and zero Neumann conditions on top and bottom boundaries. The hydraulic conductivity is a random space function with mean $\langle K\rangle=15$~m/day, with Gaussian correlation of the $\ln K$ field, correlation length $\lambda=1$~m, and variance $\sigma^2=0.1$, generated by summing 10 random modes with the Kraichnan algorithm.  An ensemble of velocity fields corresponding to 100 realizations of the $K$ field is obtained with the flow solver used in Section~\ref{2Dreg_scale_heterogeneous}, for the resolution of the GRW lattice defined by space steps $\Delta x=\Delta y=0.1$.

\begin{figure}[h]
\begin{minipage}[t]{0.45\linewidth}\centering
\includegraphics[width=\linewidth]{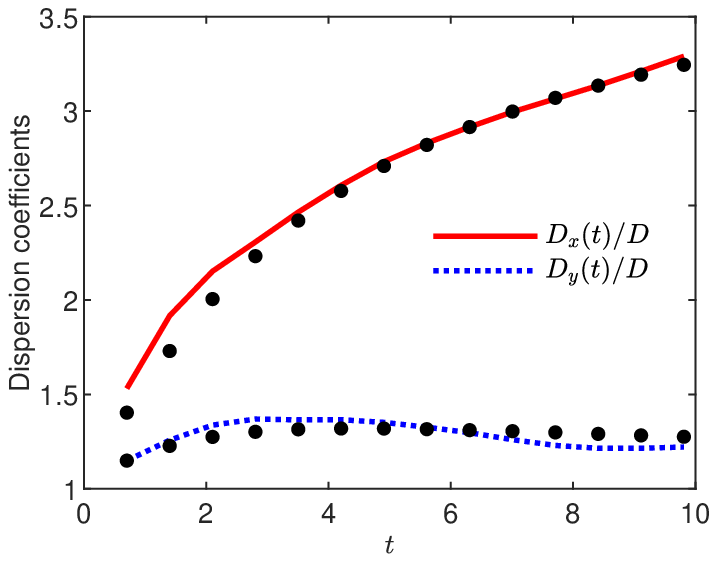}
\caption{\label{fig_DxyGRW}Dispersion coefficients estimated from GRW solutions for 100 realizations of the isotropic hydraulic conductivity $K$, with Gaussian correlated $\ln K$ field of variance $\sigma^2=0.1$ and correlation length $\lambda=1$m, in the domain $[0,20]\times[0,10]$, compared to first-order results (dots).}
\end{minipage}
\hspace*{0.1in}
\begin{minipage}[t]{0.45\linewidth}\centering
\includegraphics[width=\linewidth]{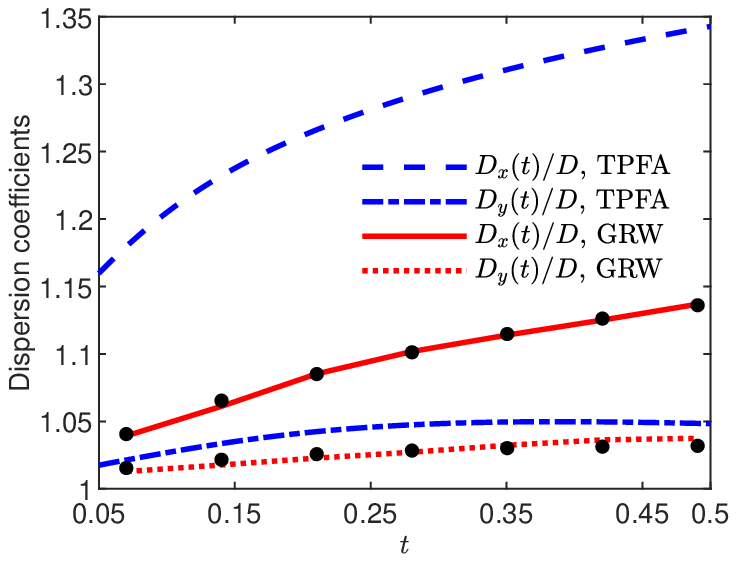}
\caption{\label{fig_compDxy}Comparison of dispersion coefficients obtained by GRW, TPFA, and first-order approximation (dots) fron an ensemble of 100 realizations of the isotropic hydraulic conductivity $K$, with Gaussian correlated $\ln K$ field of variance $\sigma^2=0.1$ and correlation length $\lambda=0.1$m, in the domain $[0,2]\times[0,1]$.}
\end{minipage}
\end{figure}
\begin{figure}[h]
\begin{minipage}[t]{0.45\linewidth}\centering
\includegraphics[width=\linewidth]{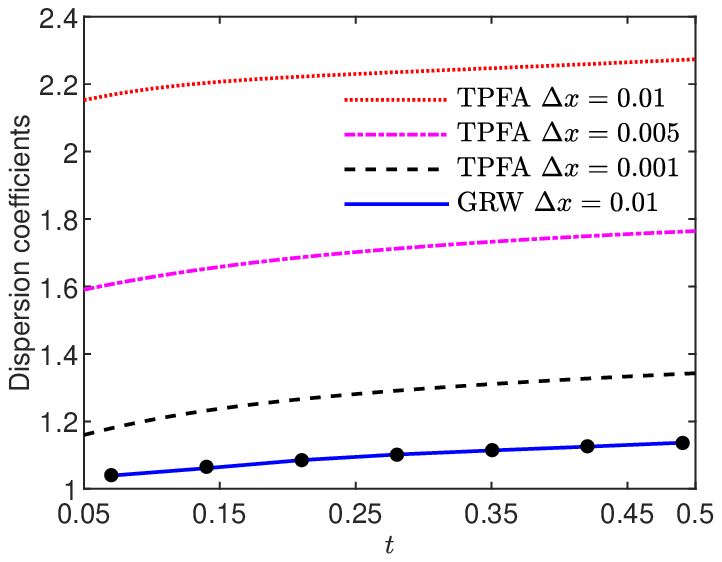}
\caption{\label{fig_DxTPFA}Comparison of longitudinal dispersion coefficients $D_x(t) / D$ obtained by GRW, TPFA, and first-order approximation (dots) from an ensemble of 100 realizations of the isotropic hydraulic conductivity $K$, with Gaussian correlated $\ln K$ field of variance $\sigma^2=0.1$ and correlation length $\lambda=0.1$m, in the domain $[0,2]\times[0,1]$.}
\end{minipage}
\hspace*{0.1in}
\begin{minipage}[t]{0.45\linewidth}\centering
\includegraphics[width=\linewidth]{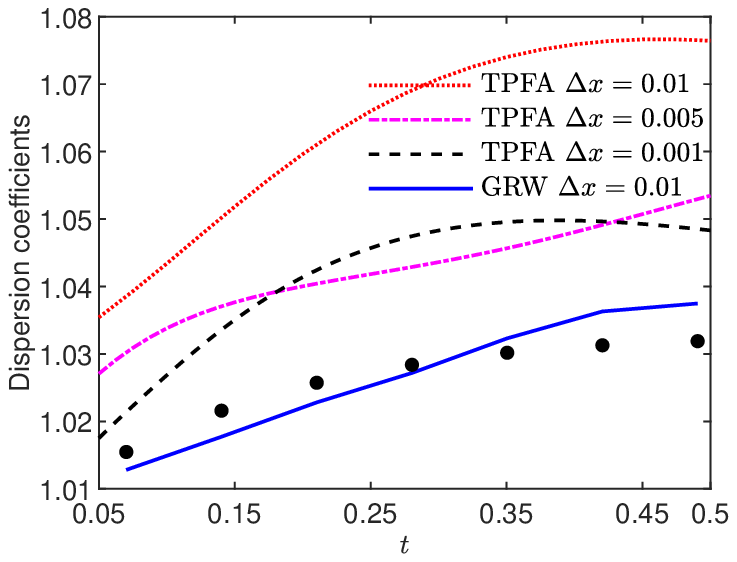}
\caption{\label{fig_DyTPFA}Comparison of transverse dispersion coefficients $D_y(t) / D$ obtained by GRW, TPFA, and first-order approximation (dots) from an ensemble of 100 realizations of the isotropic hydraulic conductivity $K$, with Gaussian correlated $\ln K$ field of variance $\sigma^2=0.1$ and correlation length $\lambda=0.1$m, in the domain $[0,2]\times[0,1]$.}
\end{minipage}
\end{figure}

Further, Monte Carlo simulations of advection-diffusion are carried out using the velocity realizations and the isotropic local dispersion coefficient $D=0.01$~m$^2$/day. The linear transport equation obtained by setting $\theta=const$ in Eq.~(\ref{eqFT02}) is solved with the unbiased GRW algorithm described in Section~\ref{2DalgFT_GRW} by using $\mathcal{N}=10^{24}$ particles to represent the concentration. The final time $T=10$ days is chosen such that the support of the concentration does not reach the boundaries during the simulation. Hence, the Monte Carlo inferences can be compared with results of linear theory which provides first-order approximations of dispersion coefficients for small variances $\sigma^2$ \cite{Bellinetal1992}. In turn, such linear approximations are accurately retrieved by averaging over ensembles of particle tracking simulations of diffusion in realizations of velocity fields approximated to the first-order in $\sigma^2$ by a Kraichnan procedure \cite{Schwarzeetal2001}. Following this approach, to infer dispersion coefficients in linear approximation, we use an ensemble of $10^4$ realizations of Krainchan velocity fields, computed with 100 random modes by the algorithm described in \cite[Appendix C.3.2.2]{Suciu2019}, and the unbiased GRW solver, with $\mathcal{N}=10^{24}$ particles in each realization. Longitudinal and transverse ``ensemble'' dispersion coefficients, $D_x$ and $D_y$, are computed as half the slope of the ensemble average of the second spatial moments  of the concentration distribution, centered at the ensemble average center of mass \cite{Bellinetal1992,Raduetal2011,Schwarzeetal2001}. The results presented in Fig.~\ref{fig_DxyGRW} show a that, in spite of relatively small ensemble of velocity realizations, the ensemble dispersion coefficients obtained with the 100 GRW solutions of the full flow problem are quite close to the reference linear results.

The computation of the velocity realizations with the transient GRW flow solver requires $10^4$ to $10^5$ iterations to fulfill the convergence criterion (\ref{eq5}) with tolerances $\varepsilon_a=\varepsilon_r=5\cdot 10^{-7}$ and about 160 seconds per realization. For the chosen discretization, $\Delta x=\Delta y=0.1$, the unbiased GRW transport solver requires, according to (\ref{eqFT2a}), a relatively rough time discretization of $\Delta t=0.5$. This leads to a total computation time  of about 1.4 seconds for the estimation of the dispersion coefficients by averaging over the 100 realizations of the statistical ensemble. By comparison, the TPFA codes needs about 3.8 seconds to compute a velocity realization and about 13 seconds for a single transport realization, by using the same spatial resolution and a time step $\Delta t=0.05$. But the TPFA estimates of the dispersion coefficients deviate by more than one order of magnitude from the linear reference solution. Since reducing the spatial steps and the local P\'{e} to reduce the numerical diffusion dramatically increases the computational burden for the TPFA codes, we solved a rescaled problem. So, to preserve the mean and the spatial variability of the velocity field, we chose a smaller domain $[0,2]\times[0.1]$, correlation length of the $\ln K$ field $\lambda=0.1$, and a new Dirichlet condition, $h(0,y)=0.1$. Now, the TPFA codes require about 60 seconds to compute one flow realization and about 3 hours for a transport realization, with $\Delta x=\Delta y=0.001$ and $\Delta t=0.0005$. The computation times for the GRW codes to solve the rescaled problem by using $\Delta x=\Delta y=0.01$ and $\Delta t=0.07$ are practically unchanged. Figure~\ref{fig_compDxy} shows that the GRW estimations of the dispersion coefficients are again close to the linear approximation. Instead the TPFA estimations deviate from the linear approximation by 10\% to 20\%. The improvements of the TPFA results with the decrease of the space steps is illustrated in Figs.~\ref{fig_DxTPFA}~and~\ref{fig_DyTPFA}. The deviations shown by the TPFA coefficients in Figs.~\ref{fig_compDxy}-\ref{fig_DyTPFA} are of the same order of magnitude as the numerical diffusion estimates (Table~\ref{table_NumDiff_TPFA}), computed as in Section~\ref{NumDiff} for constant velocity corresponding to the constant conductivity $\langle K\rangle$, in case of longitudinal coefficients $D_x$ but two order of magnitude larger in case of transverse coefficients $D_y$.
\begin{table}[!ht]
\centering
\caption{Estimation of numerical diffusion for the TPFA code.}
\label{table_NumDiff_TPFA}
\begin{tabular}{c c c c c c  c c c}
\hline
   & $\Delta x$ & $\Delta t$ &  P\'{e}  & $\varepsilon_{D_{x}}$ & $\varepsilon_{D_{y}}$ \\
  \hline
   &  0.1   & 0.05    & 7.50  & 1.84e+0 & 6.47e-2 &\\
   &  0.05  & 0.025   & 3.75  & 1.03e+0 & 3.37e-2 &\\
   &  0.01  & 0.005   & 0.75  & 1.02e+0 & 7.00e-3 &\\
   &  0.005 & 0.0025  & 0.375 & 5.06e-1 & 3.50e-3 &\\
   &  0.001 & 0.0005  & 0.075 & 1.00e-1 & 7.07e-4 &\\
 \hline
\end{tabular}
\end{table}

\section{Conclusions}
\label{concl}

The GRW schemes for simulating flow in either unsaturated or saturated porous media are equivalent to finite-difference schemes, in their deterministic implementation, or for sufficiently large numbers of particles in randomized implementations. The same, in case of BGRW solver for transport problems. Instead, the unbiased GRW is a superposition of Euler schemes for It\^{o} equation \cite{Suciu2019}, which is no longer equivalent with a finite difference scheme, unless the coefficients of the transport equation are constant. In simulations of reactive transport, GRW algorithms can use huge numbers of computational particles, even as large as the number of molecules involved in reactions, allowing simple and intuitive representations of the process.

While unbiased GRW algorithms are mainly efficient in obtaining fast solutions for large-scale transport in aquifers, BGRW solvers are appropriate for
computing solutions of fully coupled flow and transport problems in soil systems with fine variation of the parameters. The algorithms are implemented as iterative $L$-schemes which linearize the Richards equation and describe the transition from unsaturated to saturated regime. The GRW/BGRW solutions are first-order accurate in time and second-order accurate in space. For saturated regimes, the flow solver becomes a transient scheme solving steady-state flows in aquifers.

Since the GRW algorithms are explicit schemes which do not need to solve systems of algebraic equations, they are simpler and, in some cases, faster than finite element/volume schemes. The GRW $L$-schemes for non-steady coupled problems for flow and transport in soils, as well as for transport simulations in saturated aquifers, are indeed much faster than the TPFA codes used as reference in this study. However, the flow solutions for saturated porous media in large domains (e.g. field or regional scale) require much larger computing time than classical numerical schemes, due to the large number of iterations needed to achieve the convergence of the transitory scheme used to compute steady-state solutions (see also a detailed analysis in \cite{Alecsaetal2020}).

The obvious advantage of the GRW schemes is that they are practically free of numerical diffusion. This is demonstrated by the results for decoupled transport presented in Table~\ref{table_relative_errorsNumDiff}. But, as shown by the discussion at the end of Section~\ref{coupledFT}, the flow solvers also can be affected by numerical diffusion, which is difficult to isolate from other errors occurring in coupled flow and transport problems. Such errors are avoided by GRW algorithms, which prevent the occurrence of the numerical diffusion by using consistent definitions of the jump probabilities as functions of the coefficients of the flow and transport equations.

\appendix
\section{Orders of convergence for GRW $L$-schemes}
\label{appA}

\begin{figure}[h]
\begin{minipage}[t]{0.45\linewidth}\centering
\includegraphics[width=\linewidth]{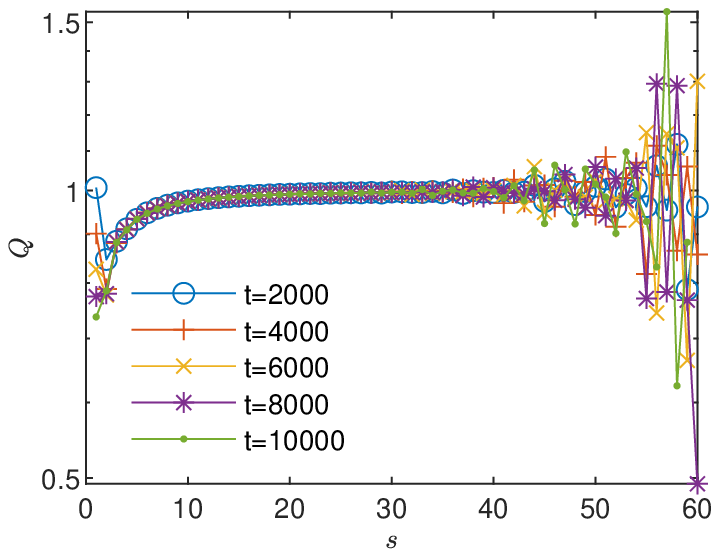}
\caption{\label{figQ_Scenario1}Computational order of convergence $Q$ (\ref{eq_Q}) of the $L$-scheme for the one-dimensional case, Scenario (1).}
\end{minipage}
\hspace*{0.1in}
\begin{minipage}[t]{0.45\linewidth}\centering
\includegraphics[width=\linewidth]{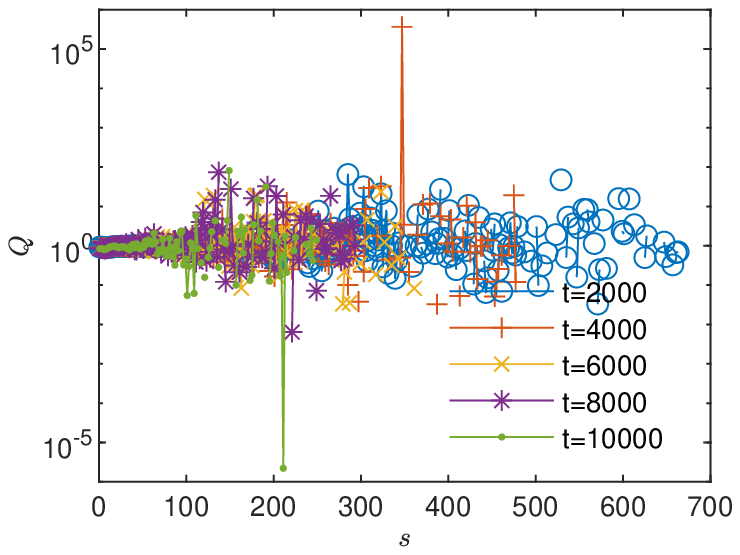}
\caption{\label{figQ_Scenario2}Computational order of convergence $Q$ (\ref{eq_Q}) of the $L$-scheme for one-dimensional case, Scenario (2).}
\end{minipage}
\end{figure}

\begin{figure}[h]
\begin{minipage}[t]{0.45\linewidth}\centering
\includegraphics[width=\linewidth]{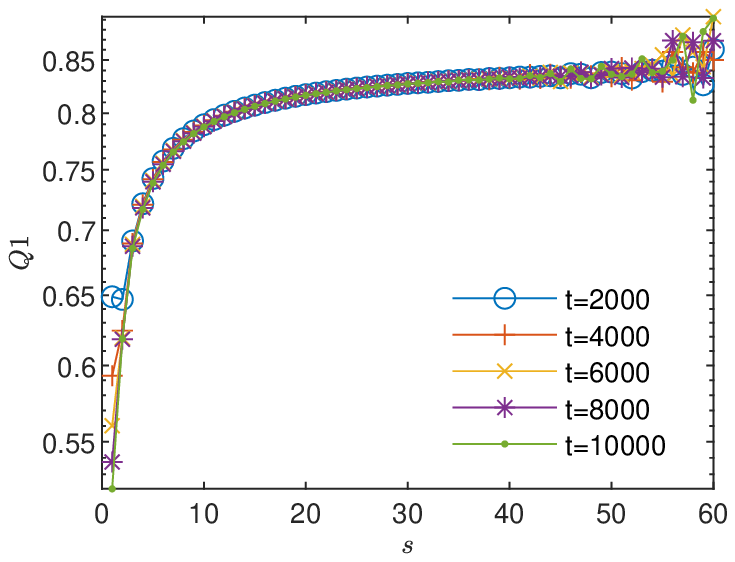}
\caption{\label{figQ1_Scenario1}Estimation of computational order of convergence $Q_1$ (\ref{eq_Q1}) of the $L$-scheme for one-dimensional case, Scenario (1).}
\end{minipage}
\hspace*{0.1in}
\begin{minipage}[t]{0.45\linewidth}\centering
\includegraphics[width=\linewidth]{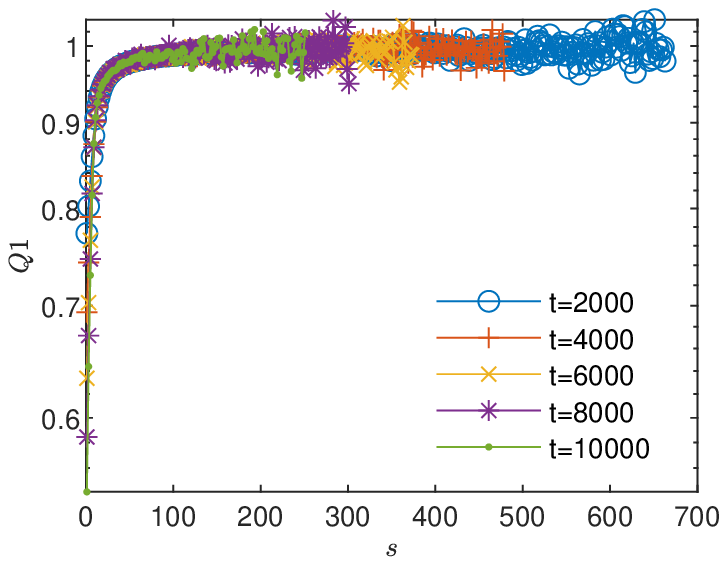}
\caption{\label{figQ1_Scenario2}Estimation of computational order of convergence $Q_1$ (\ref{eq_Q1}) of the $L$-scheme for one-dimensional case, Scenario (2).}
\end{minipage}
\end{figure}

The convergence of the iterative schemes can be investigated numerically by analyzing sequences of successive corrections, such as the norms $\|\psi^{s}-\psi^{s-1}\|$ of the convergence criterium (\ref{eq5}). The computational order of convergence $Q$ (denoted by $Q'_{\Lambda}$ in \cite{Catinas2019}) is estimated according to
\begin{equation}\label{eq_Q}
Q=\lim_{s\rightarrow\infty}\frac{\log\frac{\|\psi^{s+1}-\psi^{s}\|}{\|\psi^{s}-\psi^{s-1}\|}}
{\log\frac{\|\psi^{s}-\psi^{s-1}\|}{\|\psi^{s-1}-\psi^{s-2}\|}}.
\end{equation}
The limits (\ref{eq_Q}) for the two one-dimensional problems solved in Section~\ref{1Dvalid}, estimated in Figs.~\ref{figQ_Scenario1}~and~\ref{figQ_Scenario2} indicate the convergence of order one for both scenarios. Additionally, the estimated order of convergence is verified by employing the definition of the semicomputational order of convergence ($C'$-convergence in \cite{Catinas2020}),
\begin{equation}\label{eq_Q1}
Q_q=\lim_{s\rightarrow\infty}\frac{\|\psi^{s+1}-\psi^{s}\|}{\|\psi^{s}-\psi^{s-1}\|^q},
\end{equation}
with convergence order set to $q=1$. The convergence is called ``linear'' if the limit (\ref{eq_Q1}) verifies $Q_1<1$. Figures~\ref{figQ1_Scenario1}~and~\ref{figQ1_Scenario2} indicate the linear convergence for Scenario 1 but not for Scenario 2. In the latter case it is found that the corrections decay proportionally to $s^{-1}$, that is, the sequence converges slower than in case of linear convergence \cite[Sect. 2.2]{Catinas2020}. Results indicating linear convergence were found for both soil models in the two-dimensional flow benchmark from Section~\ref{2DvalidF}. Slower convergence with successive corrections $\sim s^{-b},\;b\in[0.1,1]$ is found in all cases of the benchmark for coupled flow and surfactant transport presented in Section~\ref{2DvalidFT}.

\section*{Acknowledgements}
The authors are grateful to Dr. Emil C\u{a}tina\c{s} for fruitful discussions on convergent sequences and successive approximation approaches. Nicolae Suciu acknowledges the financial support of the Deutsche Forschungsgemeinschaft (DFG, German Research Foundation) under Grant SU 415/4-1 -- 405338726 ``Integrated global random walk model for reactive transport in groundwater adapted to measurement spatio-temporal scales''. VISTA, a collaboration between the Norwegian Academy of Science and Letters and Equinor, funded the research of Davide Illiano, project number 6367, project name: adaptive model and solver simulation of enhanced oil recovery.


%
%
%




\end{document}